\documentclass[11pt,leqno]{article}

\usepackage{amsmath,amsfonts,amscd,amssymb,theorem}

\long\def\comment#1\endcomment{}

\comment
\endcomment

\makeatletter
\begingroup
\gdef\th@dotted{\normalfont\itshape
  \def\@begintheorem##1##2{%
        \item[\hskip\labelsep \theorem@headerfont ##1\ ##2.]}%
\def\@opargbegintheorem##1##2##3{%
   \item[\hskip\labelsep \theorem@headerfont ##1\ ##2\ (##3).]}}
\endgroup
\makeatother

\theoremstyle{dotted}

\newtheorem{theorem}{Theorem}[section]
\newtheorem{lemma}[theorem]{Lemma}

\newtheorem{prop}[theorem]{Proposition}
\newtheorem{corr}[theorem]{Corollary}

\makeatletter
\begingroup
\gdef\th@upshape{\normalfont
  \def\@begintheorem##1##2{%
        \item[\hskip\labelsep \theorem@headerfont ##1\ ##2.]}%
\def\@opargbegintheorem##1##2##3{%
   \item[\hskip\labelsep \theorem@headerfont ##1\ ##2\ (##3).]}}
\endgroup
\makeatother

\theoremstyle{upshape}

\newtheorem{defn}[theorem]{Definition}
\newtheorem{remark}[theorem]{Remark}
\newtheorem{exa}[theorem]{Example}

\makeatletter
\renewcommand{\subsection}{\@startsection{subsection}{2}{0pt}{-3ex
plus -1ex minus -0.2ex}{-2mm plus -0pt minus
-2pt}{\normalfont\bfseries}} 
\renewcommand{\subsubsection}{\@startsection{subsubsection}{3}{0pt}{-3ex
plus -1ex minus -0.2ex}{-2mm plus -0pt minus
-2pt}{\normalfont\bfseries}} 
\makeatother

\makeatletter
\@addtoreset{equation}{section}
\makeatother

\newcommand{\cntrct}                
{\hspace{2pt}\raisebox{1pt}{\text{$\lrcorner$}}\hspace{2pt}}

\newcommand{\proof}[1][Proof.]{\smallskip\noindent{\em #1}}
\def\endproof{\hfill\ensuremath{\square}\par\medskip}

\def\eqref#1{\thetag{\ref{#1}}}

\let\latexref=\ref
\def\ref#1{{\normalfont{\latexref{#1}}}}

\newcommand{\wt}{\widetilde}
\newcommand{\wh}{\widehat}

\newcommand{\idot}{{\:\raisebox{1pt}{\text{\circle*{1.5}}}}}
\newcommand{\hdot}{{\:\raisebox{3pt}{\text{\circle*{1.5}}}}}
\newcommand{\eps}{\varepsilon}
\renewcommand{\phi}{\varphi}

\newcommand{\Fun}{\operatorname{Fun}}

\newcommand{\Hom}{\operatorname{Hom}}
\newcommand{\End}{\operatorname{End}}
\newcommand{\Ext}{\operatorname{Ext}}
\newcommand{\Tor}{\operatorname{Tor}}
\newcommand{\Coker}{\operatorname{Coker}}
\newcommand{\Ker}{\operatorname{Ker}}
\renewcommand{\Im}{\operatorname{Im}}

\newcommand{\gr}{\operatorname{\sf gr}}

\newcommand{\id}{\operatorname{\sf id}}

\newcommand{\D}{\mathcal{D}}
\newcommand{\DF}{\mathcal{D}F}

\newcommand{\Z}{\mathbb{Z}}

\newcommand{\C}{\mathcal{C}}
\newcommand{\E}{\mathcal{E}}

\newcommand{\tr}{\operatorname{\sf tr}}

\newcommand{\vH}{\check{H}}
\newcommand{\vC}{\check{C}}

\newcommand{\I}{{\sf I}}

\newcommand{\ex}{\operatorname{\mathcal{E}\text{\it x}}}
\newcommand{\el}{\operatorname{\mathcal{E}\text{\it l}}}
\newcommand{\sq}{\operatorname{\mathcal{S}\text{\it q}}}
\newcommand{\spl}{\operatorname{{\mathcal{S}\text{{\it pl}\/}}}}
\newcommand{\lift}{\operatorname{\mathcal{L}\text{\it ift}}}

\newcommand{\El}{\operatorname{\text{\rm El}}}
\newcommand{\Spl}{\operatorname{\text{\rm Spl}}}
\newcommand{\Lift}{\operatorname{\text{\rm Lift}}}
\newcommand{\Sq}{\operatorname{\text{\rm Sq}}}

\newcommand{\amod}{\operatorname{\text{\rm -mod}}}
\newcommand{\bimod}{\operatorname{\text{\rm -bimod}}}
\newcommand{\proj}{\operatorname{\text{\rm -proj}}}
\newcommand{\flmod}{\operatorname{\text{\rm -flat}}}

\newcommand{\B}{\mathcal{B}}

\newcommand{\Cyl}{\operatorname{\text{\rm Cone}}}

\newcommand{\dst}{\overset{\:\raisebox{-1pt}{\text{\circle*{1.5}}}}-}

\newcommand{\motimes}{\circ}

\newcommand{\Fr}{\operatorname{\sf Fr}}

\title{Bokstein homomorphism as a universal object}

\author{D. Kaledin\thanks{Partially supported by RScF, grant
    14-21-00053, and the Dynasty Foundation award}}

\begin{document}

\maketitle

\tableofcontents

\section*{Introduction.}

For any prime $p$, the homology $H_\idot(X,k)$ of a topological
space $X$ with coefficients in the residue field $k=\Z/p\Z$ is
naturally a module over the Steenrod algebra $A$ of stable
cohomological operations at $p$. While the operations themselves are
linear endomorphisms of the functor $H_\idot(-,k)$, their
construction is highly non-linear. In particular, if one represents
$H_\idot(-,k)$ as the homology of any of the standard functorial
complexes, e.g.\ the singular chain complex $C_\idot(-,k)$, then
there is no meaningful way to lift the Steenrod algebra to a DG
algebra $A_\idot$ that acts naturally on $C_\idot(-,k)$. Indeed,
were it possible to do it, one would obtain a triangulated functor
from the stable homotopy category localized at $p$ to the derived
category $\D(A_\idot)$ of DG modules over $A_\idot$, and Adams
spectral sequence would then show that the functor is an
equivalence. Thus we would have a DG enhancement for the stable
homotopy category localized at $p$. This is known not to exist (see
e.g.\ \cite{sch} for an excelent overview, and in particular for a
proof of non-existence).

A seeming exception to the rule is the very first of the homological
operations --- the Bokstein homomorphism $\beta:H_\idot(X,k) \to
H_{\idot - 1}(X,k)$. To construct it, all one has to do is to
consider the chain complex $C_\idot(X,\Z/p^2)$ with coefficients in
the ring $\Z/p^2$. This is still a linear object, at least in the
sense of algebraic triangulated categories as in \cite{kel}, and the
Bokstein homomorphism is just the connecting homomorphism in the
long exact sequence associated to changing the coefficients from
$\Z/p^2$ to $k=\Z/p$.

However, this exception is in fact an illusion: just as all the
other homological operations, the Bokstein homomorphism does not
lift to a functorial map $C_\idot(X,k) \to C_{\idot -1}(X,k)$.

If one tries to understand what goes wrong, one realizes that the
problem is in the ring map $a:\Z/p^2 \to k$. This is a square-zero
extension of rings. If for a moment we replace $\Z/p^2$ with the
trivial square-zero extension $k[t]/t^2$, then for any complex
$V_\idot$ of flat $k[t]/t^2$-modules, the tensor product complex
$\overline{V}_\idot = V_\idot \otimes_{k[t]/t^2} k$ does have a
functorial endomorphism $\beta:\overline{V}_\idot \to
\overline{V}_{\idot -1}$. Moreover, sending $V_\idot$ to $\langle
\overline{V}_\idot,\beta \rangle$ descends to an equivalence
$$
\D(k[t]/t^2) \cong \D(k[\beta]),
$$
where $k[\beta]$ is the free associative DG algebra generated by one
element $\beta$ of degree $1$ (this is one of the simplest instances
of Koszul duality, see e.g.\ \cite{BBD}). Then by the general
homological yoga, the non-trivial extension $a:\Z/p^2 \to k$ ought to
correspond to a class in the second Hochschild cohomology group of
the ring $k$, and we might expect to have a version of Koszul
duality twisted by such a class. But the class does not exist --- in
fact, the group in question is trivial. While $a$ is a square-zero
extension of rings, it is not a square-zero extension of
$k$-algebras. So, while the derived category $\D(\Z/p^2)$ of
$Z/p^2$-modules of course has a DG enhancement, it does not have a
enhancement over $k$ --- something that would have been automatic by
Koszul duality were we to have a DG lifting of the Bokstein
homomorphism $\beta$.

\medskip

One can then pose the following somewhat philosophical question. Assume
that we are working with the category $C_\idot(k)$ of complexes of
$k$-vector spaces, with its derived category $\D(k)$ and the rest of
standard homological machinery. Assume at the same time that our
vision is limited to the world where $p=0$: we do not know that $\Z$
exists, nor that there are any primes. Can we nevertheless describe,
or even discover, the category $\D(\Z/p^2)$, purely in terms of
complexes of $k$-vector spaces and various operations on them?

\medskip

The question of course might look spurious; our justification for
posing it is that it seems to admit an interesting answer. A
one-sentence formulation is that one should simply allow
non-linearity, as long as we have linearity on the level of
homology. Specifically, we mean the following. Any DG algebra
$A_\idot$ over $k$ defines an endofunctor $F$ of the category
$C_\idot(k)$ by setting
$$
F(V_\idot) = A_\idot \otimes_k V_\idot,
$$
and since $A_\idot$ is a DG algebra, the functor $F$ is a monad
(that is, an algebra in the monoidal category of endofunctors, with
the monoiadl structure given by composition). The functor $F$ is by
definition linear, so that it sends quasiisomorphic complexes to
quasiisomorphic complexes, and induces a $k$-linear endofunctor of
the derived category $\D(k)$. However, the converse it not true:
there are endofunctors of the category $C_\idot(k)$ with these two
properties that do not come from any complex $A_\idot$. On the level
of homology, such endofunctors are necessarily linear, but on the
level of complexes, they are not. Nevertheless, they can be monads,
too; and it is such a monad that should give a functorial lifting of
the Bokstein homomorphism to chain complexes. After localization
with respect to quasiisomorphisms, the category of algebras over
such a monad would then give the derived category $\D(\Z/p^2)$.

\medskip

Of course, the above is just a sketch of a possible argument. But
note that our question can already be posed on the level of objects,
without going to complexes: how do we describe modules over $\Z/p^2$
in terms of vector spaces over $k$? More concretely, can we make
such a description functorial enough so that it allows us, for
example, to treat square-zero extension of the form $\Z/p^2 \to k$
on a par with usual $k$-linear square-zero extensions of
$k$-algebra? In this context, a lot is already known since the
pioneering work of M. Jibladze and T. Pirashvili in the late 80-ies
(see e.g.\ \cite{pira} that also contains references to earlier
works). In particular, the appropriate replacement of the second
Hochschild homology group $HH^2(R,M)$ of a $k$-algebra $R$ with
coefficients in an $R$-bimodule $M$ is the so-called {\em MacLane
  cohomology group} $H^2_M(R,M)$. By themselves, MacLane cohomology
groups have been discovered by MacLane \cite{mcl} back in 1956, but
the definition used an explicit complex. Among other things,
Jibladze and Pirashvili show that just as Hochschild cohomology
groups $HH^\hdot(R,M)$, the groups $H^\hdot_M(R,M)$ can be
interepreted as $\Ext$-groups from $R$ to $M$. To do this, one needs
to enlarge the category where the $\Ext$-groups are
computed. Specifically, $R$-bimodules correspond to additive
endofunctors of the additive category $\E=R\proj$ of finitely
generated projective $R$-modules, and instead of them, one should
use the category of all endofunctors.

Recently, an extension of this story to a more general exact or
abelian category $\E$ has been studied in \cite{ka-lo}, and several
MacLane-type cohomology theories have been constructed. A natural
next thing to do would be to show that second MacLane cohomology
classes correspond to square-zero extensions. In the case of rings
and projective modules over them, this has also been shown by
Jibladze and Pirashvili in \cite{pira}. However, their approach was
focused on rings and the associated ``algebraic theories'' rather
that just the category $\E=R\proj$, and it is not clear how to
generalize it to the abelian category case. Thus a simple purely
categorical construction of the correspondence would be useful.

\medskip

The first goal of the present paper is to provide such a
construction. We start with a ring $R$ and a bimodule $M$, interpret
MacLane cohomology classes $\eps \in H^2_M(R,M)$ as $\Ext$-classes in
the category of endofunctors of $R\proj$, and show how they
correspond to square-zero extensions of $R$ by~$M$.

In addition to this, we also consider the category $C^{pf}_\idot(R)$
of perfect complexes of $R$-modules. We introduce a convenient class
of endofunctors of $C^{pf}_\idot(R)$ that we call admissible; these
are ``additive on the level of homology''. We then show that
extensions in this category of admissible endofunctors are also
related to square-zero extensions of the ring $R$.

Our motivating example for working with complexes is a certain
natural admissible endofunctor $C_\idot$ of the category
$C_\idot^{pf}(k)$ that appeared recently in the work on cyclic
homology done in \cite{ka3}. The precise definition of the functor
$C_\idot$ is given in Subsection~\ref{cycl.dg.subs}. It is based on
an old observation: for any $k$-vector space $V$ and any integer
$i$, we have a functorial isomorphism
$$
\vH_i(\Z/p\Z,V^{\otimes p}) \cong V,
$$
where $\vH_\idot$ stands for Tate homology, and the cyclic group
$\Z/p\Z$ acts on the $p$-th tensor power $V^{\otimes p}$ by the
longest permutation. In \cite{ka3}, this isomorphisms is extended to
complexes $V_\idot \in C^{pf}_\idot(k)$, and this gives rise to the
admissible endofunctor $C_\idot$.

For applications to cyclic homology, it is important to figure out
how $C_\idot$ is related to the square-zero extension $\Z/p^2 \to
k$. This is the second goal of our paper.

What we discover along the way is that $C_\idot$ can be generalized
to an endofunctor of the category $C_\idot^{pf}(R)$ for an arbitrary
commutative ring $R$, and the corresponding square-zero extension is
then the second Witt vectors ring $W_2(R)$ (if $R=k=\Z/p$, then
$W_2(R)=\Z/p^2$). We note that while the Witt vectors ring is a very
classic object, it continues to attract attention (see for example a
recent paper \cite{den}). Thus yet another interpretation of
$W_2(R)$ might prove useful.

Finally, in the interests of full disclosure, we should mention that
our results for complexes are not as strong as for objects, and were
we to restrict our attention to objects, the paper could have been
shortened quite significantly. Nevertheless, we do spend some time
on the case of complexes. Aside from the pragmatic application to
cyclic homology and Witt vectors, we feel that it is important to
study this case in as much detail as possible, since potentially, it
provides technology that can be used for higher homological
operations. Of course, those can be also studied by a variety of
standard methods of algebraic topology, but it never hurts to have
yet another point of view. In particular, polynomial endofunctors of
the category of $k$-vector spaces have been the subject of a lot of
research in the last decades, and we now have a lot of structural
information on them (see e.g.\ \cite{poly} for an
overview). Constructing non-trivial admissible polynomial
endofunctors of $C_\idot(k)$ similar to our cyclic powers functor
$C_\idot$ and feeding this information back into algebraic topology
might be a worthwhile exercise.

\medskip

Here is an outline of the paper. Section~\ref{ele.sec} is
essentially devoted to linear algebra: we assemble several standard
facts on objects and complexes in abelian categories, in order to
set up the definitions and notation. Mostly, these deal with
representing classes in $\Ext^2$ by extensions and splitting these
extensions. All the proofs are about as difficult as the snake
lemma, so we omit them. In Section~\ref{fun.sec}, we fix an
associative unital ring $R$, and we discuss the categories of
endofunctors $\B(R)$, $\B_\idot(R)$ of the category $R\proj$ of
finitely generated projective $R$-modules and the category
$C^{pf}_\idot(R)$ of bounded complexes in $R\proj$. For $R\proj$,
simply considering pointed functors is good enough. For
$C^{pf}_\idot(R)$, we introduce the notion of an admissible
endofunctor (Definition~\ref{adm.def}). This notion is rather {\em
  ad hoc} and probably not the best possible, but it works for our
purposes, so we leave it at that. In Section~\ref{sq.sec}, we
introduce square-zero extensions of the ring $R$, and show how to
classify them by extensions in the category $\B(R)$. We also explain
why if one wants an absolute theory --- that is, square-zero
extensions of rings and not of algebras --- then it is necessary to
consider non-additive functors. The main technical point in the
correspondence is Proposition~\ref{Z.spl.prop} that shows that a
certain explicit functorial extension of abelian groups
automatically splits. It might have an abstract interpretation as a
vanishing of a certain MacLane-type cohomology group, but we have
not pursued this: instead, we opt for a direct construction. After
that, in Section~\ref{spl.sec}, we assume known the correspondence
between square-zero extension and extensions in the endofunctor
category, and we show how this allows to describe modules and
complexes over a square-zero extension $R'$ of a ring $R$ in terms
of data purely in the categories $R\proj$ and $C_\idot(R)$. Then in
Section~\ref{mult.sec}, we study multiplication: we assume that the
ring $R$ is commutative, and show that its commutative square-zero
extensions $R'$ correspond to extensions in $\B(R)$ that are
multiplicative in a certain natural way. In particular, this allows
us to show that liftings of an $R$-algebra $A$ to a square-zero
extension $R'$ are classified by Hochschild cohomology classes. We
then attempt to generalize this to complexes and DG algebras, and
explain why it does not work too well (this is
Subsection~\ref{dg.mult.subs}).

Finally, in Section~\ref{cycl.sec}, we turn to our specific example
of a non-linear extension that appears in the study of cyclic
homology: the functor $C_\idot$ corresponding to the second Witt
vectors ring $W_2(R)$. In this special case, we are able to go
further than in the general situation, and in particular, we can
describe liftings of a DG algebra $A_\idot$ over a perfect field $k$
of characteristic $p$ in terms of splittings of a certain canonical
square-zero extension of $A_\idot$ in the category of DG algebras
over $k$.

\subsection*{Acknowledgements.}

Over the years, I have benefitted from discussing the subject of the
present paper with many people. I want to mention specifically
firstly, my co-author Wendy Lowen, without whom this paper certainly
would not exist, and secondly, Vadim Vologodsky, to whom I
essentially owe the idea for Proposition~\ref{lft.ten.prop}. By
myself, I certainly would not be brave enough to imagine that a
statement of this sort could be true.

\section{Elementary extensions.}\label{ele.sec}

\subsection{Objects.}

We start by recalling some standard linear algebra. For any abelian
category $\E$, we will denote by $C_\idot(\E)$ the category of chain
complexes in $\E$ (our indexes are homological, so that the
differential acts as $d:M_{i+1} \to M_i$). For any complex $E_\idot$
and integer $m$, we denote by $E_\idot[m] \in C_\idot(\E)$ the
homological shift by $m$ of the complex $E_\idot$ (that is, the
$n$-th term of the complex $E_\idot[m]$ is $E_{n-m}$). For any
complex $E_\idot$, we denote by $\Cyl(E_\idot) \in C_\idot(\E)$ the
cone of the identity map $\id:E_\idot \to E_\idot$, so that we have
a short exact sequence
\begin{equation}\label{cone.eq}
\begin{CD}
0 @>>> E_\idot @>{\beta}>> \Cyl(E_\idot) @>{\alpha}>>
E_\idot[1] @>>> 0
\end{CD}
\end{equation}
in $C_\idot(\E)$. For any integer $m$, we will denote by $C_{\geq
  m}(\E),C_{\leq m}(\E) \subset C_\idot(\E)$ the full subcategories
spanned by complexes $E_\idot$ such that $E_i=0$ for $i < m$
resp.\ $i > m$. For any two integers $m \leq n$, we will denote
$C_{[m,n]}(\E) = C_{\geq m}(\E) \cap C_{\leq n}(\E)$, and we will
say that a complex $E_\idot$ has {\em amplitude} $[m,n]$ if it lies
in $C_{[m,n]}(\E) \subset C_\idot(\E)$. Any object $A \in \E$ can be
treated as a complex by placing it in degree $0$; this gives an
equivalence
\begin{equation}\label{e.ce}
\E \cong C_{[0,0]}(\E) \subset C_\idot(\E).
\end{equation}
For any integer $m$, we denote by $\tau_{\geq m}:C_\idot(\E) \to
C_\idot(\E)$ the composition of the embedding $C_{\geq m}(\E)
\subset C_\idot(\E)$ with its right-adjoint functor, and we denote
by $\tau_{\leq m}:C_\idot(\E) \to C_\idot(\E)$ the composition of
the embedding $C_{\leq m}(\E) \subset C_\idot(\E)$ with its
left-adjoint functor. For any two integers $m \leq n$, we denote
$\tau_{[m,n]} = \tau_{\leq n}\tau_{\geq m} \cong \tau_{\geq
  m}\tau_{\leq n}$. These are the {\em canonical truncation}
functors. For any complex $E_\idot$ and any integer $m$,
$\tau_{[m,m]}E_\idot \in C_{[m,m]}(\E) \cong \E$ is the $m$-th
homology object $H_m(E_\idot)$ of the complex $E_\idot$. We record
right away the following obvious fact.

\begin{lemma}\label{triv.le}
Assume given an exact sequence
$$
\begin{CD}
0 @>>> B_\idot @>>> C_\idot @>>> A_\idot @>>> 0
\end{CD}
$$
of complexes in $\E$, and assume that for some integer $m$, the
connecting differential $\delta:H_m(A_\idot) \to H_{m-1}(B_\idot)$
of the corresponding long exact sequence of homology is equal to
$0$. Then the sequences
$$
\begin{CD}
0 @>>> \tau_{\geq m]}B_\idot @>>> \tau_{\geq m}C_\idot @>>>
\tau_{\geq m}A_\idot @>>> 0,\\
0 @>>> \tau_{\leq m-1}B_\idot @>>> \tau_{\leq m-1}C_\idot @>>>
\tau_{\leq m-1}A_\idot @>>> 0
\end{CD}
$$
are exact.
\end{lemma}

\proof{} Clear. \endproof

\begin{defn}\label{qex.def}
A sequence
\begin{equation}\label{qex.eq}
\begin{CD}
0 @>>> B_\idot @>{b}>> E_\idot @>{a}>> A_\idot @>>> 0
\end{CD}
\end{equation}
of complexes in an abelian category $\E$ is {\em quasiexact} if $b$
is injective, $a$ is surjective, and the complex $\Ker a/\Im b$ is
acyclic.
\end{defn}

We note right away that a quasiexact sequence generates a long exact
sequence in homology: if for some integer $i$, we denote by $Z_i
\subset E_i$ the preimage $Z_i = d^{-1}(b(B_{i-1}))$ of $b(B_{i-1})
\subset E_{i-1}$, then the natural map $Z_i \to H_i(A_\idot)$
induced by $a$ is surjective, and the natural map $Z_i \to
H_{i-1}(B_\idot)$ induced by the differential $d$ factors through
$H_i(A_\idot)$.

\begin{defn}
An {\em elementary extension} of an object $A \in \E$ by an object
$B \in \E$ is a quasiiexact sequence
\begin{equation}\label{ele.ex}
\begin{CD}
0 @>>> B[1] @>{b}>> C_\idot @>{a}>> A @>>> 0
\end{CD}
\end{equation}
in $C_{[0,1]}(\E)$, where $A$ and $B$ are treated as complexes via
the embedding \eqref{e.ce}.
\end{defn}

Equivalently, an elementary extension \eqref{ele.ex} is given by two
objects $C_0$, $C_1$ in the category $\E$ and a four-term exact
sequence
\begin{equation}\label{4term}
\begin{CD}
0 @>>> B @>{b}>> C_1 @>{d}>> C_0 @>{a}>> A @>>> 0,
\end{CD}
\end{equation}
where $d:C_1 \to C_0$ is the differential in the complex
$C_\idot$. Thus by Yoneda, it defines an element $\psi \in
\Ext^2(A,B)$. Elementary extensions of $A$ by $B$ form a category
$\el(A,B)$ in an obvious way, and we denote by $\El(A,B)$ the set of
its connected components. We note that {\em a priori} it is a class,
not a set, but in all our examples it will actually be a set, and
the same will be true for various other connected component sets
that will appear. The set $\El(A,B)$ is naturally identified with
$\Ext^2(A,B)$ --- that is, elementary extensions have the same class
$\psi$ if and only if they can be connected by a chain of maps. In
this case, we say that elementary extensions are {\em equivalent}. A
canonical elementary extension with trivial class $\psi=0$ is
obtained by setting $C_1=B$, $C_0=A$, $d=0$. An elementary extension
is {\em split} if it is equivalent to the trivial one.

\begin{defn}\label{spl.def}
A {\em splitting} of an elementary extension \eqref{ele.ex} is an
object $C_{01} \in \E$ equipped with two maps $c_1:C_1 \to C_{01}$,
$c_0:C_{01} \to C_0$ that fit into a cartesian commutative diagram
$$
\begin{CD}
C_1 @>>> C_{01}\\
@VVV @VVV\\
\Im d @>>> C_1,
\end{CD}
$$
where $\Im d$ is the image of the differential $d:C_1 \to C_0$.
\end{defn}

Equivalently, a splitting is an object $C_{01}$ equipped with a
three-step filtration $W_0C_{01} \subset W_1C_{0,1} \subset C_{01}$
and identifications
$$
W_1C_{01} \cong C_1, \qquad C_{01}/W_0C_{01} \cong C_0
$$
such that under these identifications, $W_0C_{01} \subset W_1C_{01}$
is identified with $B \subset C_1$, the quotient $C_{01}/W_1C_{01} =
(C_{01}/W_0C_{01})/(W_1C_{01}/W_0C_{01})$ is identified with $B =
C_1/\Im d$, and the natural map $W_1C_{01} \to C_{01} \to
C_{01}/W_0C_{01}$ is identified with the differentil $d:C_1 \to
C_0$.  It is well-known and easy to check that an elementary
extension is split if and only if it admits a splitting in the sense
of Definition~\ref{spl.def}.

Morphisms of splittings of a fixed elementary extention $\phi =
\langle C_\idot,a,b \rangle \in \el(A,B)$ are defined in the obvious
way, and they are obviously invertible. We denote by $\spl(\phi)$
the groupoid of such splittings, and we denote by $\Spl(\phi)$ the
set of isomorphism classes of objects in $\spl(\phi)$.

The category $\E^o$ opposite to $\E$ is also abelian, and
$C_\idot(\E^o) \cong C_\idot(\E)$ is opposite to
$C_\idot(\E)$. Since the diagram \eqref{4term} is obviously
self-dual, an elementary extension $\phi \in \el(A,B)$ defines an
elementary extension of $B$ by $A$ in the opposite category, so that
we have a natural antiequivalence
\begin{equation}\label{el.du}
\el(A,B) \cong \el(B^o,A^o), \qquad \phi \mapsto \phi^o
\end{equation}
where $A^o$ and $B^o$ are $A$ and $B$ considered as objects of
$\E^o$, and $\phi^o$ is the extension given the sequence
\eqref{4term} corresponding to $\phi$ but considered as a sequence
in $\E^o$. Alternatively, in terms of the diagram \eqref{ele.ex},
one needs to pass to the opposite category and apply a homological
shift by $1$. The notion of a splitting is also self-dual, so that
we have a natural equivalence
$$
\spl(\phi) \cong \spl(\phi^o), \qquad \phi \in \el(A,B)
$$
sending $\langle C_{01},c_0,c_1 \rangle$ to $\langle
C_{01}^o,c_1,c_0 \rangle$.

Finally, assume given an elementary extension $\phi = \langle
C_\idot,a,b \rangle \in \el(A,B)$, and a map $f:B \to B'$ to some
other object $B' \in \E$. Then we define the composition $f \circ
C_\idot$ by the short exact sequence
$$
\begin{CD}
0 @>>> B[1] @>{-f \oplus b}>> B'[1] \oplus C_\idot @>>> f \circ
C_\idot @>>> 0
\end{CD}
$$
in $C_\idot(\E)$, and we note that $f \circ C_\idot$ is naturally an
elementary extension of $A$ by $B'$. The construction is functorial,
so that we have a functor
\begin{equation}\label{comp.l.ex}
f \circ -:\el(A,B) \to \el(A,B').
\end{equation}
On the level of the connected component sets, this functor
corresponds to the map $f \circ -:\Ext^2(A,B) \to \Ext^2(A,B')$
obtained by composition with $f$. Moreover, sending a splitting
$\langle C_{01},c_0,c_1 \rangle$ of the elementary extension $\phi$
to the cokernel of the natural map $(-f)\oplus (c_1 \circ b):B \to
B' \oplus C_{01}$ defines a functor
\begin{equation}\label{comp.l.spl}
f \circ -:\spl(\phi) \to \spl(f \circ \phi).
\end{equation}
Dually, for any map $f:A' \to A$, one obtains a functor
\begin{equation}\label{comp.r.ex}
- \circ f:\el(A,B) \to \el(A',B)
\end{equation}
that induces the composition map $- \circ f:\Ext^2(A,B) \to
\Ext^2(A',B)$ on the connected component sets, and one has a natural
functor
\begin{equation}\label{comp.r.spl}
- \circ f:\spl(\phi) \to \spl(\phi \circ f)
\end{equation}
for any elementary extension $\phi \in \el(A,B)$.

\subsection{Complexes.}

More generally, assume given two complexes $A_\idot,B_\idot \in
C_\idot(\E)$, and define an {\em elementary extension of $A_\idot$
  by $B_\idot$} as a quasiexact sequence
\begin{equation}\label{ele.ex.c}
\begin{CD}
0 @>>> B_\idot[1] @>{b}>> C_\idot @>{a}>> A_\idot @>>> 0
\end{CD}
\end{equation}
in $C_\idot(\E)$. Elementary extensions again form a category in the
obvious way; we denote this category by
$\el_\idot(A_\idot,B_\idot)$, and we denote by $\El_\idot(A,B)$ the
set of its connected components. Extensions are {\em equivalent} if
they lie in the same connected component of
$\el_\idot(A_\idot,B_\idot)$. A trivial extension is obtained by
taking $C_\idot = A_\idot \oplus B_\idot$. An extension is {\em
  split} if it is equivalent to the trivial one.

Every elementary extension \eqref{ele.ex.c} defines an extension of
$A_\idot$ by the complex $\Ker a \subset C_\idot$, and $\Ker a$ is
by definition quasiisomorphic to $B_\idot[1]$, so that we again
obtain a canonical element
$$
\psi \in \Hom_{\D(\E)}(A_\idot,B_\idot[2]),
$$
where $\D(\E)$ denotes the derived category of the abelian category
$\E$. Sending an extension \eqref{ele.ex.c} to its class $\psi$
identifies the set $\El_\idot(A_\idot,B_\idot)$ with the group
$\Hom_{\D(\E)}(A_\idot,B_\idot[2])$.

We note that if $A_\idot$ and $B_\idot$ are objects $A$, $B$ placed
in degree $0$, then $\el_\idot(A,B)$ is bigger than $\el(A,B)$ since
the complex $C_\idot$ in \eqref{ele.ex.c} is not required to have
amplitude $[0,1]$. However, we have a fully faithful embedding
$\el(A,B) \subset \el_\idot(A,B)$, and replacing a complex $C_\idot$
with its canonical truncation $\tau_{[0,1]}C_\idot$ defines a
natural functor
$$
\tau_{[0,1]}:\el_\idot(A,B) \to \el(A,B)
$$
inverse to the embedding. Both the embedding and the inverse
functor induce an isomorphism on the sets of connected components,
so that every elementary extension in $\El_\idot(A,B)$ is
canonically equivalent to an extension of the form
\eqref{ele.ex}. The notion of an elementary extension of complexes
is also self-dual, so that we have a natural antiequivalence
$$
\el_\idot(A_\idot,B_\idot)^o \cong \el_\idot(A_\idot^o,B_\idot^o),
$$
an extension of \eqref{el.du}.

\begin{defn}\label{dg.spl.def}
A {\em left DG splitting} of an elementary extension $\phi$
represented by a diagram \eqref{ele.ex.c} is a complex $C^l_\idot$
equipped with maps $l:C^l_\idot \to C_\idot$, $b^l:\Cyl(B_\idot) \to
C^l_\idot$ so that, with the map $\alpha$ being the map of
\eqref{cone.eq}, the diagram
\begin{equation}\label{dg.spl.sq}
\begin{CD}
\Cyl(B_\idot) @>{b^l}>> C_\idot^l\\
@V{\alpha}VV @VV{l}V\\
B_\idot[1] @>{b}>> C_\idot
\end{CD}
\end{equation}
is commutative, and the sequence
\begin{equation}\label{dg.spl.seq}
\begin{CD}
0 @>>> \Cyl(B_\idot) @>{b^l}>> C^l_\idot @>{a \circ l}>> A_\idot
@>>> 0
\end{CD}
\end{equation}
is quasiexact. A {\em right DG splitting} of the elementary
extension $\phi$ is a complex $C^r_\idot$ equipped with maps
$r:C_\idot \to C_\idot^r$, $a^r:C_\idot^r \to \Cyl(A_\idot)$ such
that $\langle C^r_\idot,r,a^r\rangle$ give a left DG splitting of
the extension $\phi^o$ in the opposite category $C_\idot(\E)^o$.
\end{defn}

As a corollary of the definition, we see that for left DG splitting
$C^l_\idot$, the map $a \circ l:C^l_\idot \to A_\idot$ is a
quasiisomorphism, and for a right DG splitting $C^r_\idot$, the map
$r \circ b:B_\idot[1] \to C^r_\idot$ is a quasiisomorphism. An
elementary extension admits a right DG splitting iff it admits a
left DG splitting iff it is split. For any left DG splitting
$\langle C^l_\idot,l,b^l \rangle$ of an elementary extension, the
cone $C^{lr}_\idot$ of map $l$ is a right DG splitting of the same
extension, with the obvious maps $r$, $a^r$. Conversely, for a right
DG splitting $\langle C^r_\idot,r \rangle$, the shift $C^{rl}_\idot$
of the cone of the map $r$ is a left DG splitting. If we denote by
$\spl_\idot^l(\phi)$, $\spl^r_\idot(\phi)$ the categories of left
and right DG splittings of an elementary extension $\phi$, then this
defines functors
\begin{equation}\label{lr.eq}
\spl^l_\idot(\phi) \to \spl^r_\idot(\phi), \qquad
\spl^r_\idot(\phi) \to \spl^l_\idot(\phi).
\end{equation}
The functors are not mutually inverse, but we do have functorial
quasiisomorphisms
$$
C^{rlr}_\idot \cong C^r_\idot, \qquad C^{lrl}_\idot \cong C^l_\idot.
$$
Thus if we denote by $\Spl_\idot^l(\phi)$
resp.\ $\Spl_\idot^r(\phi)$ the sets of connected components of the
categories $\spl_\idot^l(\phi)$ resp.\ $\spl_\idot^r(\phi)$, then we
have a natural identification $\Spl_\idot^l(\phi) \cong
\Spl_\idot^r(\phi)$.

\begin{defn}
A left DG splitting $\langle C_\idot^l,l,b^l \rangle$ is {\em
  strict} if the diagram \eqref{dg.spl.sq} is a Cartesian square in
$C_\idot(\E)$. A right DG splitting is {\em strict} if it is strict
as a left DG splitting in $C_\idot(\E^o)$.
\end{defn}

For a strict left DG splitting $\langle C^l_\idot,l,b^l \rangle$ of
an extension \eqref{ele.ex.c}, the middle homology of the quasiexact
sequence \eqref{dg.spl.seq} is canonically isomorphic to the middle
homology of the sequence \eqref{ele.ex.c}. For a general left DG
splitting, we just have a map between two acyclic complexes, and not
necessarily an isomorphism. The same is true for right DG
splittings. Any map between strict DG splittings is automatically an
isomorphism, so that, if we denote by
$$
\spl^l(\phi) \subset \spl^l_\idot(\phi), \qquad \spl^r(\phi) \subset
\spl^r_\idot(\phi)
$$
the full subcategories spanned by strict splittings, then both
$\spl^l(\phi)$ and $\spl^r(\phi)$ are groupoids.

For any two objects $A,B \in \E$, a splitting $\langle
C_{01},c_0,c_1 \rangle$ of an elementary extension \eqref{ele.ex} in
the sense of Definition~\ref{spl.def} defines its left and right DG
splittings $C^l_\idot$, $C^r_\idot$ in the sense of
Definition~\ref{dg.spl.def}: these are the natural complexes
\begin{equation}\label{dg.lr}
\begin{CD}
C_1 @>{c_1}>> C_{01},
\end{CD}
\qquad
\begin{CD}
C_{01} @>{c_0}>> C_0
\end{CD}
\end{equation}
placed in homological degrees $0$ and $1$, with the maps $l$
resp.\ $r$ induced by the maps $c_0$, $\id$ resp.\ $\id$, $c_1$, and
the maps $b^l$, $a^l$ induced by the maps $b$, $c_1 \circ b$
resp.\ $a \circ c_1$, $a$. Both these DG splittings are strict and
functorial, so that we obtain natural functors
$$
\spl(\phi) \to \spl^l(\phi) \subset \spl^l_\idot(\phi), \qquad
\spl(\phi) \to \spl^r(\phi) \subset \spl^r_\idot(\phi).
$$
Automatically, $C^r_\idot$ is quasisomorphic to $C^{lr}_\idot$, and
$C^l_\idot$ is quasiisomorphic to $C^{rl}_\idot$, but both these
quasiisomorphisms are not isomorphisms on the nose.

Conversely, given a left DG splitting $\langle C^l_\idot,l,b^l
\rangle$ of an elementary extension \eqref{ele.ex}, we can
canonically construct its splitting $C_{01}$ in two steps. First, we
observe that the canonical truncation $\wt{C}^l_\idot =
\tau_{[0,1]}C^l_\idot$ is also a left DG splitting, with the maps
$\wt{l} = \tau_{[0,1]}(l)$, $\wt{b}^l = \tau_{[0,1]}(b^l)$. Second,
we define $C_{01}$ by the short exact sequence
\begin{equation}\label{DG.el.spl}
\begin{CD}
0 @>>> \wt{C}^l_1 @>{-d \oplus \wt{l}}>> \wt{C}^l_0 \oplus C_1 @>>>
C_{01} @>>> 0,
\end{CD}
\end{equation}
where $d:\wt{C}^l_1 \to \wt{C}^l_0$ is the differential. The map
$c_1$ is induced by the natural embedding $C_1 \to \wt{C}^l_0 \oplus
C_1$, and the map $c_0$ is induced by the map
$$
\begin{CD}
\wt{C}^l_0 \oplus C_1 @>{\wt{l} \oplus d}>> C_0,
\end{CD}
$$
where $d$ is the differential in the complex $C_\idot$ (since $l$
commutes with differentials, this map factors through $C_{01}$). The
original left DG splitting $C^l_\idot$ is then automatically
equivalent, but not necessarily isomorphic, to the left-hand side
complex in \eqref{dg.lr}. For right DG splittings, one has a
parallel construction that we leave to the reader.

\subsection{Gerbs.}

For any objects $A,B \in \E$ and an elementary extension $\phi \in
\el(A,B)$, the set $\Spl(\phi)$ of isomorphism classes of splittings
of $\phi$ is naturally a torsor over the extension group
$\Ext^1(A,B)$. This can be lifted to the level of categories: if we
denote by $\ex(A,B)$ the groupoid of short exact sequences
\begin{equation}\label{ext.dia}
\begin{CD}
0 @>>> B @>{b}>> E @>{a}>> A @>>> 0
\end{CD}
\end{equation}
in the category $\E$, then $\spl(\phi)$ is a gerb over
$\ex(A,B)$. Specifically, for any two objects $\langle
C_{01}',c_1'c_0' \rangle, \langle C_{01}'',c_1'',c_0''\rangle \in
\Spl(\phi)$, the diagram
$$
\begin{CD}
C_1 @>{c_1' \oplus c_1''}>> C'_{01} \oplus C_{01}'' @>{c_0' \oplus
  -c_0''}>> C_0
\end{CD}
$$
is a complex, and its middle homology is naturally an extension of
$A$ by $B$, so that we obtain a difference functor
\begin{equation}\label{diff.fu}
-:\spl(\phi) \times \spl(\phi) \to \ex(A,B).
\end{equation}
Analogously, for any $\langle C_{01},c_1,c_0 \rangle \in
\spl(\phi)$ and any short exact sequence \eqref{ext.dia}, the
diagram
$$
\begin{CD}
B @>{c_1 \oplus b}>> C_{01} \oplus E @>{c_0 \oplus -a}>> A
\end{CD}
$$
is a complex, and its middle homology is naturally an object in
$\spl(\phi)$, so that we obtain a difference functor
\begin{equation}\label{sum.fu}
-:\spl(\phi) \times \ex(A,B) \to \spl(\phi).
\end{equation}
The standard algebraic relations between differences then extend to
isomorphisms between the relevant compositions of the functors
\eqref{diff.fu} and \eqref{sum.fu}. In particular, we have natural
isomorphisms
\begin{equation}\label{sum.diff}
C \cong C' - (C' - C), \qquad E \cong C - (C - E),
\end{equation}
for any $C,C' \in \spl(\phi)$, $E \in \ex(A,B)$.

To extend the functors \eqref{diff.fu}, \eqref{sum.fu} to DG
splittings, for any two complexes $A_\idot,B_\idot \in C_\idot(\E)$,
we denote by $\ex_\idot(A_\idot,B_\idot)$ the category of quasiexact
sequences \eqref{qex.eq}. Then for any left DG splitting $\langle
C^l_\idot,l,b^l\rangle$ of a DG elementary extension $\phi = \langle
C_\idot,a,b \rangle \in \el_\idot(A_\idot,B_\idot)$, and any object
$\langle E_\idot,a',b' \rangle \in \ex_\idot(A_\idot,B_\idot)$, we
have a natural sequence
$$
\begin{CD}
0 @>>> B_\idot @>{(b^l \circ \beta) \oplus b'}>> C^l_\idot \oplus
E_\idot @>{(a \circ l) \oplus -a'}>> A_\idot @>>> 0
\end{CD}
$$
that is exact on the left and on the right, and its middle homology
is naturally a left DG splitting of the extension $\phi$. We denote
this splitting by $C^l_\idot - E_\idot$. The construction is
functorial, so that we obtain a difference functor
\begin{equation}\label{dg.sum.fu}
-:\spl^l_\idot(\phi) \times \ex_\idot(A_\idot,B_\idot) \to
\spl^l_\idot(\phi).
\end{equation}
We note that since any exact sequence of complexes is quasiexact, we
have a natural fully faithful embedding $\ex(A_\idot,B_\idot)
\subset \ex_\idot(A_\idot,B_\idot)$, and if a DG splitting $C$ is
strict, then for any $E \in \ex(A_\idot,B_\idot) \subset
\ex_\idot(A_\idot,B_\idot)$, the difference $C - E$ is also strict,
so that \eqref{dg.sum.fu} induces a functor
\begin{equation}\label{dg.sum.fu.1}
-:\spl^l(\phi) \times \ex(A_\idot,B_\idot) \to
\spl^l(\phi).
\end{equation}
Conversely, assume given another left DG splitting $\langle
\wt{C}^l_\idot,\wt{l},\wt{b}^l\rangle \in \spl^l_\idot(\phi)$, and
denote by $\overline{C}_\idot$ the cone of the map
\begin{equation}\label{c.wc.cone}
\begin{CD}
C^l_\idot \oplus \wt{C}^l_\idot @>{l \oplus -\wt{l}}>> C_\idot.
\end{CD}
\end{equation}
Then by definition, we have a quasiexact sequence
\begin{equation}\label{ol.qex}
\begin{CD}
0 @>>> \overline{B}_\idot @>>> \overline{C}_\idot @>>>
\overline{A}_\idot @>>> 0,
\end{CD}
\end{equation}
where $\overline{B}_\idot$ is the cone of the map
$$
\begin{CD}
\Cyl(B_\idot) \oplus \Cyl(B_\idot) @>{\alpha \oplus -\alpha}>>
B_\idot[1],
\end{CD}
$$
and $\overline{A}_\idot$ is the cone of the map
$$
\begin{CD}
A_\idot \oplus A_\idot @>{\id \oplus -\id}>> A_\idot.
\end{CD}
$$
Furthermore, $\id \oplus \id:\Cyl(B_\idot) \to \Cyl(B_\idot) \oplus
\Cyl(B_\idot)$ extends to a natural embedding $\Cyl(B_\idot) \subset
\overline{B}_\idot[-1]$, and the projection $\id \oplus -\id:A_\idot
\oplus A_\idot \to A_\idot$ extends to a natural surjection
$\overline{A}_\idot \to \Cyl(A_\idot)$. Thus altogether, the
sequence \eqref{ol.qex} induces a sequence
$$
\begin{CD}
0 @>>> \Cyl(B_\idot) @>>> \overline{C}_\idot[-1] @>>>
\Cyl(A_\idot)[-1] @>>> 0.
\end{CD}
$$
Its middle homology is then naturally an object in
$\ex_\idot(A_\idot,B_\idot)$ that we denote by $C^l_\idot -
\wt{C}^l_\idot$. Again, the construction is functorial, so that we
obtain a natural functor
\begin{equation}\label{dg.diff.fu}
-:\spl^l_\idot(\phi) \times \spl^l_\idot(\phi) \to
\ex_\idot(A_\idot,B_\idot).
\end{equation}
If DG splittings $C^l_\idot$, $\wt{C}^l_\idot$ are strict, then one
can refine the construction by using the kernel of the map
\eqref{c.wc.cone} instead of its cone. What we end up with, then, is
the middle cohomology of a natural sequence
$$
\begin{CD}
0 @>>> \Cyl(B_\idot) @>{b^l \oplus \wt{b}^l}>> C^l_\idot \oplus
\wt{C}^l_\idot @>{l \oplus -\wt{l}}>> C_\idot @>>> 0.
\end{CD}
$$
We denote it by $C_\idot^l \dst \wt{C}^l_\idot$. It fits naturally
into a quasiexact sequence
$$
\begin{CD}
0 @>>> B_\idot @>>> (C_\idot^l \dst \wt{C}^l_\idot) @>>> A_\idot
@>>> 0
\end{CD}
$$
whose middle homology is naturally identified with the middle
homology of \eqref{ele.ex.c} homologically shifted by $1$. In
particular, $C^l_\idot \dst \wt{C}^l_\idot$ is naturally an object
in $\ex_\idot(A_\idot,B_\idot)$. Moreover, we have a natural map
\begin{equation}\label{df.ddf}
C^l_\idot \dst \wt{C}^l_\idot \to C^l_\idot - \wt{C}^l_\idot 
\end{equation}
induced by the embedding of the kernel of a surjective map of
complexes into the homlogical shift of its cone.

The functors \eqref{dg.sum.fu} and \eqref{dg.diff.fu} are related by
natural maps
\begin{equation}\label{dg.sum.diff}
C^l_\idot \to \wt{C}^l_\idot - (\wt{C}^l_\idot - C^l_\idot), \qquad
E_\idot \to C^l_\idot - (C^l_\idot - E_\idot),
\end{equation}
a DG refinement of the isomorphisms \eqref{sum.diff}. All of this
has a completely parallel version for right DG splittings that we
leave to the reader.

\section{Functors.}\label{fun.sec}

\subsection{Pointed functors.}

For any small category $\C$ and any abelian category $\E$, the
category $\Fun(\C,\E)$ of all functors from $\C$ to $\E$ is
abelian. If the category $\C$ is additive, then say that a functor
$F \in Fun(\C,\E)$ is {\em pointed} if $F(0)=0$, and denote by
\begin{equation}\label{point.eq}
\Fun_o(\C,\E) \subset \Fun(\C,\E)
\end{equation}
the full subcategory spanned by pointed functors. Note that since $0
\in \C$ is both the initial and the terminal object, every $F
\in \Fun(\C,\E)$ canonically splits as
\begin{equation}\label{point.F}
F = F' \oplus F(0),
\end{equation}
where $F'$ is pointed and $F(0)$ is the constant functor. Thus
$\Fun_o(\C,\E) \subset \Fun(\C,\E)$ is an abelian subcategory closed
under extensions. Another remark is that for any two objects
$V_1,V_2 \in \C$, both $V_1$ and $V_2$ are distiguished in $V_1
\oplus V_2$ as the images of orthogonal idempotents $e_1,e_2 \in
\End(V_1 \oplus V_2)$, $e_1^2=e_1$, $e_2^2=e_2$, $e_1e_2=e_2e_1=0$,
and for any pointed functor $F:\C \to \E$, the idempotents $F(e_1)$
and $F(e_2)$ are also orthogonal. Therefore we have a functorial
decomposition
\begin{equation}\label{cross.eq}
F(V_1 \oplus V_2) = F(V_1) \oplus F(V_2) \oplus F(V_1,V_2),
\end{equation}
where $F(-,-):\C \times \C \to \E$ is a certain functor known as
{\em cross-effect functor}. A pointed functor $F$ is additive if and
only if the cross-effects functor $F(-,-)$ is trivial; therefore the
full subcategory
\begin{equation}\label{add.o}
\Fun_{add}(\C,\E) \subset \Fun_o(\C,\E)
\end{equation}
is again an abelian subcategory closed under extensions.

Assume given an associative unital ring $R$. Denote by $R\amod$ the
abelian category of left $R$-modules, and let $R\proj \subset
R\amod$ be the full subcategory spanned by finitely generated
projective modules. Then the category $R\proj$ is small and
additive, and the category $R\amod$ is abelian, so we can take
$\C=R\proj$ and $\E=R\amod$. We denote
\begin{equation}\label{b.r}
\B(R) = \Fun_o(R\proj,R\amod).
\end{equation}
More generally, given another associative unital ring $R'$, we
denote
\begin{equation}\label{b.rr}
\B(R,R') = \Fun_o(R\proj,R'\amod) \subset \Fun(R\proj,R'\amod).
\end{equation}
Any left module $M$ over the product $R' \otimes R^o$ of $R'$ and
the algebra $R^o$ opposite to $R$ defines a functor $\I(M) \in
\B(R,R')$ by setting
\begin{equation}\label{i.m}
\I(M)(V) = M \otimes_R V, \qquad V \in R\proj.
\end{equation}
The functor $\I(M)$ is additive. Conversely, for any additive
functor $F:R\proj \to R'\amod$, the value $F(R)$ at the free
$R$-module $R$ is naturally a module over $R' \otimes R^o$, and we
have $F \cong \I(F(R))$. Thus $(R' \otimes R^o)\amod \cong
\Fun_{add}(R\proj,R'\amod)$, and \eqref{add.o} induces an exact
fully faithful embedding
$$
\I(-):(R' \otimes R^o)\amod \subset \B(R,R')
$$
whose essential image is closed under extensions. If $R=R'$, what we
get is a functor
\begin{equation}\label{i.dash}
\I(-):R\bimod \subset \B(R),
\end{equation}
where $R\bimod$ is the category of $R$-bimodules. This is a fully
faithful additive exact functor between abelian categories, and its
essential image consists of additive functors $F:R\proj \to
R\amod$. In addition, this essential image is closed under
extensions, so that the natural map
$$
\Ext^i(M,N) \to \Ext^i(\I(M),\I(N)), \qquad M,N \in R\bimod
$$
is an isomorphism for $i=0$ and $i=1$. For $i \geq 2$, this not
necessarily true.

To simplify notation, we denote $\I=\I(R) \in \B(R)$, where $R$ is
the diagonal bimodule. The object $\I$ corresponds to the natural
embedding functor $R\proj \to R\amod$.

\begin{defn}\label{ext.R.def}
For any associative unital ring $R$, and any $R$-bimodule $M$, an
{\em elementary extension} of $R$ by $M$ is an elementary extension
\eqref{ele.ex} of $\I$ by $\I(M)$ in the abelian category
\eqref{b.r}.
\end{defn}

By definition, an elementary extension $\phi$ of
Definition~\ref{ext.R.def} has a cohomology class lying in the group
\begin{equation}\label{H.M}
H^2_M(R,M) = \Ext^2_{\B(R)}(\I,\I(M)).
\end{equation}
This group coincides with the so-called second {\em MacLane
  cohomology group} of $R$ with coefficients in $M$. Indeed, by
virtue of the canonical decomposition \eqref{point.F}, the embedding
\eqref{point.F} induces a fully faithful functor on the derived
categories, so that we can compute the right-hand side of
\eqref{H.M} in the abelian category $\Fun(R\proj,R\amod)$; the
result coincides with MacLane cohomology by the famous result of
Jibladze and Pirashvili \cite{pira}.

\subsection{Admissible functors.}

In order to obtain a version of Definition~\ref{ext.R.def} for
complexes, we need to introduce an appropriate class of functors
between them. For any associative unital ring $R$, we denote by
$C_\idot(R) = C_\idot(R\amod)$ the category of complexes of left
$R$-modules, and let $C_\idot^{pf}(R) \subset C_\idot(R)$ be the
full subcategory spanned by finite-length complexes of finitely
generated projective modules. Inverting quasiisomorphisms in
$C_\idot(R)$ gives the derived category $\D(R)$ of the abelian
category $R\amod$, and $C^{pf}_\idot(R) \subset C_\idot(R)$ spans
the full subcategory $\D^{pf}(R) \subset \D(R)$ of compact objects
in the triangulated category $\D(R)$.

For any rings $R$, $R'$, the category $C^{pf}_\idot(R)$ is additive
and small, and the category $C_\idot(R')$ is abelian, so that we can
consider pointed functors from $C^{pf}_\idot(R)$ to
$C_\idot(R')$. Note that we have
$$
\Fun_o(C_\idot^{pf}(R),C_\idot(R')) \cong
C_\idot(\Fun_o(C_\idot^{pf}(R),R'\amod)).
$$
We are interested in functors that moreover descend
to the derived categories, so we need to impose some further
conditions. The following seems a reasonable choice.

\begin{defn}\label{adm.def}
A pointed functor $F_\idot:C_\idot^{pf}(R) \to C_\idot(R')$ is {\em
  admissible} if it sends acyclic complexes in $C^{pf}_\idot(R)$ to
acyclic complexes in $C_\idot(R')$, and exact sequences of complexes
in $C^{pf}_\idot(R)$ to quasiexact sequences of complexes in
$C_\idot(R')$. The category
\begin{equation}\label{b.d.rr}
\B_\idot(R,R') \subset \Fun_o(C_\idot^{pf}(R),C_\idot(R'))
\cong C_\idot(\Fun_o(C_\idot^{pf}(R),R'\amod))
\end{equation}
is the full subcategory spanned by admissible functors.
\end{defn}

Admissibility is not invariant under quasiisomorphisms. However, the
subcategory $\B_\idot(R,R') \subset
C_\idot(\Fun_o(C_\idot^{pf}(R),R'\amod))$ is obviously closed under
extensions and homological shifts, thus also under taking cones of
maps. Therefore the full subcategory
\begin{equation}\label{d.adm}
\D^{adm}(R,R') \subset \D(\Fun_o(C_\idot^{pf}(R),R'\amod))
\end{equation}
spanned by admissible functors is triangulated.

\begin{lemma}\label{adm.le}
Any admissible functor $F_\idot:C_\idot^{pf}(R) \to
C_\idot(R')$ sends quasiisomorphisms to quasiisomorphisms, hence
descends to a functor
\begin{equation}\label{adm.desc}
F:\D^{pf}(R) \to \D(R'),
\end{equation}
and the functor $F$ is triangulated.
\end{lemma}

\proof{} For any map $f:E'_\idot \to E_\idot$ in $C^{pf}_\idot(R)$,
the exact sequence
\begin{equation}\label{cone.seq}
\begin{CD}
0 @>>> E'_\idot @>>> \Cyl(f) @>>> E_\idot[1] @>>> 0
\end{CD}
\end{equation}
in $C^{pf}_\idot(R)$ maps to the exact sequence \eqref{cone.eq}, so
that we have a commutative diagram
\begin{equation}\label{dia.b}
\begin{CD}
0 @>>> F_\idot(E'_\idot) @>>> F_\idot(\Cyl(f)) @>>>
F_\idot(E_\idot[1]) @>>> 0\\
@. @V{F_\idot(f)}VV @VVV @|\\
0 @>>> F_\idot(E_\idot) @>{F_\idot(\beta)}>> F_\idot(\Cyl(E_\idot))
@>{F_\idot(\alpha)}>> F_\idot(E_\idot[1])
@>>> 0
\end{CD}
\end{equation}
in $C_\idot(R')$ with quasiexact rows. If $f$ is a quasiisomorphism,
then the middle terms are acyclic, so that $F_\idot(f)$ is indeed a
quasiisomorphism. Hence $F_\idot$ descends to a functor
\eqref{adm.desc}. Moreover, even without $E'_\idot$ and the map $f$,
we still have the second row of \eqref{dia.b}, and it induces a
diagram
$$
\begin{CD}
F_\idot(E_\idot)[1] @<<< \Cyl(F_\idot(\beta)) @>>>
F_\idot(E_\idot[1])
\end{CD}
$$
in $C_\idot(R')$. The diagram is functorial in $E_\idot$, and both
maps are quasiisomorphisms. Therefore the functor $F$ of
\eqref{adm.desc} commutes with homological shifts. Furthermore, for
any $E_\idot,E'_\idot \in C_\idot^{pf}(R)$, the image under
$F_\idot$ of the split exact sequence
\begin{equation}\label{sum.seq}
\begin{CD}
0 @>>> E_\idot @>>> E_\idot \oplus E'_\idot @>>> E'_\idot @>>> 0
\end{CD}
\end{equation}
must be quasiexact, so that the cross-effects complex
$F_\idot(E_\idot,E'_\idot)$ is acyclic, and $F$ is
additive. Finally, since every distinguished triangle in
$\D^{pf}(R)$ can be represented by an exact sequence
\eqref{cone.seq}, the diagram \eqref{dia.b} shows that $F$ sends
distinguished triangles to distinguished triangles.
\endproof

By virtue of Lemma~\ref{adm.le}, the category $\D^{adm}(R,R')$ of
\eqref{d.adm} acts naturally by triangulated functors from
$\D^{pf}(R)$ to $\D(R')$. In this sense, it can serve as a
reasonable axiomatization of the ``triangulated category of enhanced
triangulated functors from $\D^{pf}(R)$ to $\D(R')$''.

\subsection{Expansion and restriction.}

One source of examples of admissible functors is the following. For
any complex $F_\idot$ of pointed functors from $R\proj$ to
$R'\amod$, define its {\em expansion}
$$
\eps_\idot(F_\idot):C_\idot^{pf}(R) \to C_\idot(R')
$$
by applying $F_\idot$ termwise and taking the sum-total complex of
the resulting bicomplex --- that is, we set
\begin{equation}\label{e.dot}
\eps_\idot(F_\idot)(E_\idot)_n = \bigoplus_{i+j=n}F_i(E_j), \qquad E_\idot
\in C_\idot(R).
\end{equation}
The differential in the bicomplex squares to $0$ precisely because
$F_i$ are pointed functors.

\begin{defn}\label{qadd.def}
A complex $F_\idot$ in $\B(R,R')$ is {\em quasiadditive} if for any
$V,V' \in R\proj$, the complex $F_\idot(V,V')$ of cross-effects
functors is acyclic.
\end{defn}

\begin{lemma}\label{add.adm.le}
Assume given a quasiadditive complex $F_\idot$ in $\B(R,R')$. Then
its expansion $F'_\idot = \eps_\idot(F_\idot)$ is admissible in the
sense of Definition~\ref{adm.def}.
\end{lemma}

\proof{} Since by definition, complexes in $C^{pf}_\idot(R)$ are
made out of projective $R$-modules, every exact sequence in
$C^{pf}_\idot(R)$ is termwise-split. Then to show that $F'_\idot$
sends it to a quasiexact sequence, it suffices to apply the standard
spectral sequence of the bicomplex defining $\eps_\idot(F_\idot)$.

It remains to show that $F'_\idot$ sends acyclic complexes to
acyclic complexes. But every acyclic complex in $C^{pf}_\idot(R)$ is
a sum of complexes of the form $\Cyl(E)[i]$ for some integer $i$ and
projective $R$-module $E \in R\proj$. Since $F'_\idot$ obviously
commutes with shifts, it suffices to prove that $F'_\idot(\Cyl(E))$
is acyclic. This is also obvious, since by definition, we have
$F'_\idot(\Cyl(E)) \cong \Cyl(F'_\idot(E))$.
\endproof

Since \eqref{e.dot} is functorial with respect to $F_\idot$,
Lemma~\ref{add.adm.le} provides a natural functor
\begin{equation}\label{e.c.b}
\eps_\idot:C_\idot^{add}(\B(R,R')) \to \B_\idot(R,R'),
\end{equation}
where $C_\idot^{add}(\B(R,R')) \subset C_\idot(\B(R,R'))$ denotes
the full subcategory spanned by quasiadditive
complexes. Quasiadditivity is by definition invariant under
quasiisomorphisms, so that inverting quasiisomorphisms in
$C_\idot^{add}(\B(R,R'))$, we obtain a full triangulated subcategory
$$
\D^{add}(R,R') \subset \D(\B(R,R'))
$$
in the derived category $\D(\B(R,R'))$. Expansion then descends to a
functor
$$
\eps:\D^{add}(R,R') \to \D^{adm}(R,R'),
$$
where $\D^{adm}(R,R')$ is the triangulated category
\eqref{d.adm}. It seems likely that the functor $\eps$ is actually
an equivalence --- in particular, up to a quasiisomorphism, every
admissible functor from $C_\idot^{pf}(R)$ to $C_\idot(R')$ comes
from a quasiadditive complex of pointed functors from $R\proj$ to
$R'\amod$. However, we have not pursued this (in any case, the
construction of such a complex is unlikely to be either explicit or
particularly useful).

On the chain level, the expansion functor \eqref{e.c.b} is
definitely not essentially surjective. In particular, for any
pointed functor $F_\idot:C_\idot^{pf}(R) \to C_\idot(R')$ and any
integer $i$, define the {\em homological twist} $\tau^i(F_\idot)$ by
\begin{equation}\label{tw.eq}
\tau^i(F_\idot)(V_\idot) = F_\idot(V_\idot[i])[-i].
\end{equation}
Then if $F_\idot$ is admissible, $\tau^i(F_\idot)$ is also
admissible and induces the same functor from $\D^{pf}(R)$ to
$\D(R')$. However, $F_\idot$ and $\tau^i(F_\idot)$ are in general
different as functors from $C_\idot^{pf}(R)$ to $C_\idot(R')$. But
if $F_\idot = \eps_\idot(F'_\idot)$ for some quasiadditive complex
$F'_\idot$, they are the same --- as mentioned in the proof of
Lemma~\ref{add.adm.le}, we have a natural isomorphism
\begin{equation}\label{e.tw}
\tau^i(\eps_\idot(F'_\idot)) \cong \eps_\idot(F'_\idot)
\end{equation}
for any complex $F'_\idot \in C_\idot^{add}(\B(R,R'))$.

\medskip

We can also go in the other direction: composing a pointed functor
$F_\idot:C_\idot^{pf}(R) \to C_\idot(R')$ with the embedding
\eqref{e.ce} gives an object
$$
\rho_\idot(F_\idot) \in \Fun_o(R\proj,C_\idot(R')) \cong
C_\idot(\B(R,R')),
$$
and it is easy to see from \eqref{sum.seq} that if $F_\idot$ is
admissible, $\rho_\idot(F_\idot)$ is quasiadditive. We thus obtain a
{\em restriction functor}
\begin{equation}\label{p.c.b}
r_\idot:\B_\idot(R,R') \to C_\idot^{adm}(\B(R,R')).
\end{equation}
This is a one-sided inverse to the expansion functor \eqref{e.c.b}
--- we have a natural identification $\eps_\idot \circ \rho_\idot
\cong \id$. We cannot say anything about the composition in the
other direct. In particular, for an arbitrary $F_\idot \in
\B_\idot(R,R')$, we do not even have a natural map between $F_\idot$
and $\eps_\idot(\rho_\idot(F_\idot))$.

\subsection{DG elementary extensions.}

We can now give our DG version of Definition~\ref{ext.R.def}. For
any associative unital ring $R$ and $R$-module $M$, the functor
$\I(M) \in \B(R)$ is additive, hence quasiadditive in the sense of
Definition~\ref{qadd.def}. We simplify notation by letting
$\B_\idot(R) = \B_\idot(R,R)$, and we denote by
$$
\I_\idot(M) = \eps_\idot(\I(M)) \in \B_\idot(R)
$$
the expansion \eqref{e.c.b} of the functor $\I(M)$. If $M$ is the
diagonal bimodule $R$, we further simplify notation by writing
$\I_\idot = \I_\idot(R) = \eps_\idot(\I)$.

\begin{defn}\label{DG.ext.R.def}
For any associative unital ring $R$, and any $R$-bimoudle $M$, an
{\em elementary DG extension} of $R$ by $M$ is a quasiexact sequence
\begin{equation}\label{DG.R.eq}
\begin{CD}
0 @>>> \I_\idot(M)[1] @>>> C_\idot @>>> \I_\idot @>>> 0
\end{CD}
\end{equation}
in the category $\B_\idot(R) \subset
C_\idot(\Fun_o(C_\idot^{pf}(R),R\amod))$ such that for any integer
$m$, the functor $C_\idot:C_\idot^{pf}(R) \to C_\idot(R)$ sends the
subcategory $C_{\leq m}^{pf}(R) \subset C_\idot^{pf}(R)$ into the
subcategory $C_{\leq m+1}(R) \subset C_\idot(R)$.
\end{defn}

In other words, an elementary DG extension is an elementary
extension \eqref{ele.ex.c} of $\I_\idot$ by $\I_\idot(M)$ with
admissible $C_\idot$ satisfying the additional assumption: we
require that
\begin{equation}\label{bnd.1}
C_\idot(C_{\leq m}^{pf}(R)) \subset C_{\leq m+1}(R)
\end{equation}
for any integer $m$. Note that since \eqref{DG.R.eq} is required to
be quasiexact, and $\I_\idot$, $\I_\idot(M)$ are admissible, one
does not have to check all the conditions of
Definition~\ref{adm.def} to see that the functor
$C_\idot:C^{pf}_\idot(R) \to C_\idot(R)$ is admissible. In fact, it
suffices to check that it sends termwise-split injections to
termwise-split injections, and termwise-spit surjections to
termwise-split surjections; the rest is automatic.

In particular, if $C_\idot$ in \eqref{DG.R.eq} is an expansion of
some complex in the category $\Fun_o(C_\idot^{pf}(R),R\amod)$, it is
automatically admissible. Therefore any elementary extension of $R$
by $M$ in the sense of Definition~\ref{ext.R.def} generates an
elementary DG extension given by
$$
\begin{CD}
0 @>>> \I_\idot(M) = \eps_\idot(\I(M)) @>>> \eps_\idot(C_\idot) @>>>
\I_\idot = \eps_\idot(\I) @>>> 0,
\end{CD}
$$
where $\eps_\idot$ is the expansion functor
\eqref{e.c.b}. Conversely, every elementary DG extension
\eqref{DG.R.eq} gives an elementary extension
$$
\begin{CD}
0 @>>> \I(M) @>>> \tau_{\geq 0}\rho_\idot(C_\idot) @>>> \I @>>> 0,
\end{CD}
$$
where $\rho_\idot$ is the restriction functor \eqref{p.c.b}, and
$\tau_{\geq 0}$ is the canonical truncation at $0$ (the complex
$\rho_\idot(C_\idot)$ automatically lies in $C_{\leq 1}(\B(R))$ by
\eqref{bnd.1}, and the truncation is needed to insure that it also
lies in $C_{\geq 0}(\B(R))$. Both elementary extensions and
elementary DG extensions form categories in the obvious way, denoted
$\el(R,M)$ resp.\ $\el_\idot(R,M)$, and the correspondence between
them is sufficiently functorial to give expansion and restriction
functors
\begin{equation}\label{e.r.ext}
\eps:\el(R,M) \to \el_\idot(R,M), \qquad \rho:\el_\idot(R,M) \to \el(R,M).
\end{equation}
As for the functors \eqref{e.c.b} and \eqref{p.c.b}, we have $\rho
\circ \eps \cong \id$. Slightly more generally, for any integer $i$,
the homological twist $\tau^i(C_\idot)$ of an elementary DG
extension $C_\idot$ is also an elementary DG extension, and we can
define a twisted restriction functor $\rho_i:\el_\idot(R,M) \to
\el(R,M)$ by
\begin{equation}\label{tw.ext}
\rho_i(C_\idot) = \rho(\tau^i(C_\idot)).
\end{equation}
Then by \eqref{e.tw}, we still have $\rho_i \circ \eps \cong
\id$. We denote the sets of connected components of the categories
$\el(R,M)$, resp.\ $\el_\idot(R,M)$ by $\El(R,M)$,
resp.\ $\El_\idot(R,M)$.

\section{Square-zero extensions.}\label{sq.sec}

\subsection{The algebra case.}

Assume given an associative unital ring $R$ and an $R$-bimodule $M$. As
usual, by a {\em square-zero extension} $R'$ of $R$ by $M$ we will
understand an associative unital ring $R'$ equipped with a short
exact sequence
\begin{equation}\label{sq.eq}
\begin{CD}
0 @>>> M @>{i}>> R' @>{q}>> R @>>> 0
\end{CD}
\end{equation}
of abelian groups such that $q$ is a unital ring map, and we have
$$
i(m)r' = i(mp(r')), \qquad r'i(m)=i(p(r')m), \qquad i(m)i(m')=0
$$
for any $m,m' \in M$ and $r' \in R'$. Square-zero extensions of $R$
by $M$ form a category $\sq(R,M)$ in the obvious way. This category
is obviously a groupoid, and we denote by $\Sq(R,M)$ the set of
isomorphism classes of its objects.

If $R$ and $R'$ are algebras over a commutative ring $k$, and $R$ is
flat over $k$, then it is well-known that square-zero extensions are
classified by Hochschild cohomology classes: we have a natural
identification
\begin{equation}\label{sq.hh}
\Sq(R,M) \cong HH^2_k(R,M) = \Ext^2_{R \otimes_k R^o}(R,M),
\end{equation}
where the $\Ext$-group on the right is computed in the category of
$k$-linear $R$-bimodules.

Let us recall what is probably the most direct way to construct the
correspondence \eqref{sq.hh} (it seems to be a well-known folklore
result). We begin with the following general fact. For any
square-zero extension \eqref{sq.eq}, restriction with respect to $p$
gives a natural exact functor $q_*:R\amod \to R'\amod$. It has a
left-adjoint functor $q^*:R'\amod \to R\amod$. By adjunction, $q^*$
is right-exact, so it has derived functors $L^iq^*:R'\amod \to
R\amod$.

\begin{lemma}\label{L.p.le}
Assume given a square-zero extension \eqref{sq.eq} of an associative
unital ring $R$ by an $R$-bimodule $M$. Then for any projective left
$R$-module $V$, we have natural identifications
$$
q^*q_*V \cong V, \qquad L^1q^*q_*V \cong M \otimes_R V,
$$
and both are functorial in $V$.
\end{lemma}

\proof{} The first isomorphism is clear: for any $V \in R\amod$, we
have
$$
q^*q_*V \cong R \otimes_{R'} V,
$$
and the right-hand side is by definition identified with $V$. For
the second isomorphism, consider the sequence \eqref{sq.eq} as a
short exact sequence of $R'$-bimodules. Then for any $V' \in
R'\amod$, we have a long exact sequence
$$
\begin{CD}
0 = \Tor^{R'}_1(R',V') @>>> \Tor^{R'}_1(R,V') @>{\delta}>> M \otimes_{R'}
V' @>{a}>> \\
@>{a}>> R' \otimes_{R'} V' = V' @>{b}>> R \otimes_{R'} V' @>>> 0
\end{CD}
$$
functorial in $V'$. If $V' = q_*V$ for some $V \in R\amod$, then the
map $b$ is an isomorphism. Therefore the map $a$ vanishes, and the
map $\delta$ is an isomorphism. Since $L^1q^*(-) \cong
\Tor^{R'}_1(R,-)$ and $M \otimes_{R'} V \cong M \otimes_R V$, this
proves the claim.
\endproof

Assume now that the ring $R$ is a flat algebra over a commutative
ring $k$, and that a square-zero extension \eqref{sq.eq} is
$k$-linear. Then in particular, $R$ is naturally a left module over
$R' \otimes_k R^o$. Moreover, for any projective $R' \otimes_k
R^o$-module $P$ and any $R$-module $V$, the tensor product $P
\otimes_R V$ is a projective $R'$-module. Thus if we take a
projective resolution $P'_\idot$ of the $R' \otimes_k R^o$-module
$R$, then for any projective $R$-module $V$, the complex $P'_\idot(V)
= P_\idot \otimes_R V$ is a projective resolution of the $R'$-module
$q_*V$. We then obtain a functor
\begin{equation}\label{p.res}
P'_\idot(-):R\proj \to C_\idot(R'\amod)
\end{equation}
that sends any $V \in R\proj$ to a projective resolution of
$q_*V$. Composing it with $q^*$, we obtain a functor
$$
P_\idot(-):R\proj \to C_\idot(R\amod), \quad
P_\idot(V) = q^*P'_\idot(V) = R \otimes_{R'} P'_\idot \otimes_R V
$$
represented by a complex
$$
P_\idot = R \otimes_{R'} P'_\idot
$$
of modules of $R \otimes_k R^o$. Moreover, by definition, for any $V
\in R\proj$, the $i$-th homology of the complex $P_\idot(V)$ is
identified with $L^iq^*q_*V$. Thus if we let
$$
\overline{P}_\idot = \tau_{[0,1]}P_\idot
$$
be the canonical truncation of the complex $P_\idot$, then
Lemma~\ref{L.p.le} shows that $\overline{P}_\idot$ fits into en
exact sequence
$$
\begin{CD}
0 @>>> M @>>> \overline{P}_1 @>>> \overline{P}_0 @>>> R @>>> 0
\end{CD}
$$
of modules over $R \otimes_k R^o$. This exact sequence represents by
Yoneda a class in $HH^2_k(R,M)$. This is what corresponds to the
square-zero extension $R'$ under the identification \eqref{sq.hh}.

\subsection{The absolute case.}

Unfortunately, the construction presented above does not work in the
absolute case, that is, in the absence of a good base ring $k$. What
breaks down is \eqref{p.res}. In general, there is no way to lift
the functor $q_*$ to an additive functor from $R\proj$ to
$C_\idot(R'\amod)$ that sends a finitely generated projective
$R$-module $V$ to a projective resolution of $q_*V$.

\begin{exa}\label{w.2.exa}
Assume that $R = k$ is a finite field of some
characteristic $p$, and let $R' = W_2(k)$ be its ring of
second Witt vectors. Then $R\proj$ is the category of
finite-dimensional $k$-vector spaces, so that its objects are free
modules $k^n$, $n \geq 1$. For any such module $V=k^n$, we of course
have many projective resolutions over $W_2(k)$ --- for example, we
can take the minimal resolution
$$
\begin{CD}
\cdots @>{p}>> W_2(k)^n @>{p}>> W_2(k)^n @>{p}>> W_2(k)^n @>{p}>>
W_2(k)^n,
\end{CD}
$$
with all the differentials given by multiplication by $p$. However,
in order to make this functorial in $V$, one would have, at the very
least, to find a group map $GL_n(k) \to GL_n(W_2(k))$ splitting the
reduction map $GL_n(W_2(k)) \to GL_n(k)$. This is not possible
already for $n=2$ and $k = \Z/p\Z$. More generally, it is easy to
see that in fact {\em any} additive functor from $k\proj$ to
$W_2(k)\proj$, hence also to $C_\idot(W_2(k))$ is trivial
(such a functor must correspond to a module $M$ over $W_2(k) \otimes
k$ that is projective over $W_2(k)$, and since $W_2(k) \otimes k
\cong k$, this implies $M=0$).
\end{exa}

One way to correct the situation is to use non-additive functors and
MacLane cohomology instead of Hochschild cohomology. Namely, while
the minimal resolution of Example~\ref{w.2.exa} is hopeless, it is
nevertheless perfectly possible to lift the functor $\pi_*$ to a
functor \eqref{p.res} as long as we allow non-additive
functors. This is based on the following observation.

\begin{lemma}\label{ev.le}
For any object $c \in \C$ of a small category $\C$, and any
projective object $P \in \Fun(\C,\E)$ in the abelian category of
functors from $\C$ to an abelian category $\E$ satisfying $AB4^*$,
the value $P(c) \in \E$ of the functor $P$ at the object $c$ is a
projective object in $\E$.
\end{lemma}

\proof{} Define a functor $i_{c*}:\E \to \Fun(\C,\E)$ by setting
$$
i_{c*}(E)(c') = \prod_{f:c \to c'}E, \qquad E \in \E, c' \in \C,
$$
where the right-hand side stand for the product of copies of $E$
numbered by maps from $c$ to $c'$. Then $i_{c*}$ is obviously
right-adjoint to the evaluation functor $E \mapsto E(c)$, and since
$\E$ satisfies $AB4^*$, $i_{c*}$ is exact.
\endproof

By virtue of Lemma~\ref{ev.le}, we can now do the following. Assume
that $R'$ is a square-zero extension \eqref{sq.eq} of a ring $R$ by
a bimodule $M$, and consider the category $\B(R,R')$ of
\eqref{b.rr}. This is an abelian category that satisfies $AB4$ and
$AB4^*$ and has enough projectives. Composing a functor $F \in
\B(R,R')$ with $q^*$ gives a pointed functor $q^*F:R\proj \to
R\amod$, so that we have a natural functor
$$
q^*:\B(R,R') \to \B(R),
$$
and analogously, composition with $q_*$ gives a natural functor
$$
q_*:\B(R) \to \B(R,R').
$$
With this notation, the object $q_*(\I) \in \B(R,R')$ corresponds to
the functor $q_*$ itself. Choose a projective resolution $P'_\idot$
of $q_*(\I)$, consider the complex
$$
P_\idot = q^*P'_\idot
$$
in the category $\B(R)$, and let $\overline{P}_\idot =
\tau_{[0,1]}P_\idot$. Then by Lemma~\ref{ev.le}, for any $V \in
R\proj$, $P'_\idot(V)$ is a projective resolution of $q_*(V)$ in the
category $R'\amod$, so that as before, Lemma~\ref{L.p.le} shows that
we have an elementary extension
$$
\begin{CD}
0 @>>> \I(M)[1] @>>> \overline{P}_\idot @>>> \I @>>> 0
\end{CD}
$$
of $R$ by $M$ in the sense of Definition~\ref{ext.R.def}. This
extension depends both on $R'$ and on the choice of a resolution
$P'_\idot$, so it is not quite functorial in $R'$ and does not
define a functor $\sq(R,M) \to \el(R,M)$. However, since any two
projective resolutions of the same objects are quasiisomorphic, we
do have a well-defined map
\begin{equation}\label{sq.el}
\Sq(R,M) \to \El(R,M) \cong HH^2_M(R,M)
\end{equation}
on the sets of connected components. This associates a canonical
MacLane cohomology class of degree $2$ to any square-zero extension
of $R$ by $M$.

\subsection{Canonical splitting.}

We now want to prove that the map \eqref{sq.el} is in fact an
isomorphism --- that is, square-zero extensions \eqref{sq.eq}
correspond bijectively to their MacLane cohomology classes. To do
this, we first need to consider functors from $R\proj$ to abelian
groups.

Let $e:\Z \to R$ be the tautological map sending $1$ to $1$, so that
the restriction functor $e_*$ is just the forgetful functor from
$R$-modules to abelian groups. Then $e_*$ is pointed, thus gives an
object in $\B(R,\Z)$. To keep notation consistent, we extend $e_*$
to a functor
\begin{equation}\label{phi.eq}
e_*:\B(R) \to \B(R,\Z)
\end{equation}
by applying it pointwise; then the object in $\B(R,\Z)$
corresponding to $e_*$ is $e_*(\I)$, where $\I \in \B(R)$
corresponds to the taugological functor. Moreover, $e_*$ is
additive, so that $e_*(\I)$ corresponds to an $R^o$-module ---
namely, to $R$ considered as a right module over itself. Since $R$
is projective as a module over $R^o$, we have
\begin{equation}\label{ext.1.triv}
\Ext^1_{\B(R,\Z)}(e_*(\I),\I(W)) \cong \Ext^1_{R^o}(R,W) = 0
\end{equation}
for any $R^o$-module $W$ and the corresponding additive functor
$\I(W)$ in the category $\B(R,\Z)$. Somewhat surprisingly, an
analogous statement also holds for $\Ext^2(e_*(\I),-)$. In fact, we
even have a stronger statement.

\begin{prop}\label{Z.spl.prop}
Assume given a unital associative ring $R$, an $R^o$-module $W$, and
an elementary extension $\phi \in \el(e_*(\I),\I(W))$ represented by
a quasi\-exact sequence
\begin{equation}\label{phi.p.m.seq}
\begin{CD}
0 @>>> \I(W)[1] @>{a}>> C_\idot @>{b}>> e_*(\I) @>>> 0
\end{CD}
\end{equation}
in $C_\idot(\B(R,\Z))$. Then $\phi$ admits a splitting $C_{01}$,
and this splitting is unique up to an isomorphism.
\end{prop}

\proof{} Uniqueness immediately follows from
\eqref{ext.1.triv}. Indeed, for any two splittings $C_{01},C'_{01}$
of the extension $\phi$, we have
$$
C'_{01} \cong C_{01} - T
$$
for some object $T \in \ex(e_*(\I),\I(W))$, where in the right-hand
side we take the difference functor \eqref{sum.fu}. But by
\eqref{ext.1.triv}, we have $T \cong 0$ and $C'_{01} \cong C_{01}$.

To prove existence, we need to define a splitting $C_{01}(V)$ of the
elementary extension $\phi(V)$ of abelian groups for any $V \in
R\proj$, and we need to do it in a way that is functorial in
$V$. Evaluating \eqref{phi.p.m.seq} at the free module $R \in
R\proj$, we obtain in particular a surjective map $b:C_0(R) \to
e_*(R) = R$. Choose an element $s \in C_0(R)$ such that
$b(s)=1$. Since $R \in R\proj$ represents the forgetful functor
$e_*$, by Yoneda, the element $s \in C_0(R)$ defines a map
\begin{equation}\label{s.yo}
\wt{s}:V \to C_0(V), \qquad V \in R\proj
\end{equation}
of sets functorial in $V$, and since $b(s)=1$, we have $b \circ
\wt{s} \cong \id$. Let
\begin{equation}\label{p.times}
C_{01}(V) = C_1(V) \times V
\end{equation}
as a set, and define $c_1:C_1(V) \to C_{01}(V)$, $c_0:C_{01}(V) \to
C_0(V)$ by
\begin{equation}\label{c.0.1}
c_1(c) = c \times 0, \qquad c_0(c \times v) = \delta(c) + \wt{s}(v),
\qquad c \in C_1(V), v \in V,
\end{equation}
where $\delta:C_1(V) \to C_0(V)$ is the differential. To turn
$C_{01}(V)$ into the required splitting, we need to endow it with a
structure of an abelian group in such a way that $c_0$ and $c_1$ become
group maps. It suffices to construct a functorial cocycle map
\begin{equation}\label{coc}
c(-,-):V \times V \to C_1(V)
\end{equation}
such that
\begin{gather}
c(v_1,0) = c(0,v_1) = 0,\label{coc.eq.1}\\
c(v_1,v_2) = c(v_2,v_1),\label{coc.eq.2}\\
c(v_1,v_2) + v_3 + c(v_1+v_2,v_3) = c(v_1,v_2+v_3) + v_1 +
c(v_2,v_3),\label{coc.eq.3}\\
\delta(c(v_1,v_2)) = \wt{s}(v_1) +
\wt{s}(v_2) - \wt{s}(v_1+v_2)\label{coc.eq.4}
\end{gather}
for any $v_1,v_2,v_2 \in V$. Since the functor sending $V$ to the
set $V \times V$ is represented by $R^2$, it suffice to construct
the element $c(v_1,v_2) \in C_1(V)$ in the universal case $V = R^2$,
$v_1$ and $v_2$ the generators of the two copies of $R$. Then
\eqref{coc.eq.1} is equivalent to saying that $c(v_1,v_2)$ lies in
the cross-effect component $C_1(R,R) \subset C_1(R^2)$ of the
decomposition~\ref{cross.eq}. But since both $e_*(\I)$ and $\I(W)$
are additive functors, their cross-effect functors are trivial, and
since \eqref{phi.p.m.seq} is quasiexact, its cross-effect component
$$
\begin{CD}
C_1(R,R) @>{\delta}>> C_0(R,R)
\end{CD}
$$
is an acyclic complex. Therefore \eqref{coc.eq.4} uniquely
determines $c(v_1,v_2)$. Then \eqref{coc.eq.2} immediately follows
from uniqueness. To check \eqref{coc.eq.3}, it suffices to consider
the case $V=R^3$ spanned by $v_1$, $v_2$, $v_3$, and again, the
claim immediately follows from \eqref{coc.eq.4} and the acyclicity
of cross-effects.
\endproof

\subsection{Regular endomorphisms.}

We note that the splitting of Proposition~\ref{Z.spl.prop} is only
canonical in a weak sense: isomorphisms between different splittings
are not canonical. This is in fact an inherent feature of the
construction, and it is this feature that allows us to recover a
non-trivial square-zero extension from an elementary extension.

Namely, assume given an associative unital ring $R$, an $R$-bimodule
$M$, and an elementary extension $\phi \in \el(R,M)$ represented by
a quasiexact sequence
\begin{equation}\label{p.m.seq}
\begin{CD}
0 @>>> \I(M)[1] @>>> C_\idot @>>> \I @>>> 0
\end{CD}
\end{equation}
in the category $C_{[0,1]}(\B(R))$. Then applying the functor
\eqref{phi.eq} to \eqref{p.m.seq} gives an elementary extension
$$
\begin{CD}
0 @>>> e_*(\I(M)) @>>> e_*(C_\idot) @>>> e_*(\I) @>>> 0
\end{CD}
$$
in the category $\B(R,\Z)$. Denote by $C_{01}$ the canonical
splitting of this extension provided by
Proposition~\ref{Z.spl.prop}, and consider the algebra
$\End(C_{01})$ of additive endomorphisms of the functor $C_{01} \in
\B(R,\Z)$. Say that an endomorphism $r' \in \End(C_{01})$ is {\em
  regular} if it fits into a commutative diagram
\begin{equation}\label{r.r.r}
\begin{CD}
e_*(C_1) @>{c_1}>> C_{01} @>{c_0}>> e_*(C_0)\\
@V{r}VV @VV{r'}V @VV{r}V\\
e_*(C_1) @>{c_1}>> C_{01} @>{c_0}>> e_*(C_0)
\end{CD}
\end{equation}
for some element $r \in R$. The set of all regular endomorphisms
obviously forms a unital subalgebra in $\End(C_{01})$ that we denote
by $R'$.

\begin{lemma}\label{regu.le}
In the assumptions above, the algebra $R'$ is a square-zero extension
of $R$ by $M$.
\end{lemma}

\proof{} Note that since $c_0$ in \eqref{r.r.r} is surjective, $r$ is
unique. Sending an element $r' \in R'$ to this unique $r$ defines an
algebra map
\begin{equation}\label{p.rr}
q:R' \to R.
\end{equation}
The kernel of the map $q$ is the space of map $r'$ that vanish on
$\Im a$, thus factor through the quotient $e_*(\I) =
C_{01}/c_1(e_*(C_1))$, and take values in $e_*(\I(M)) = \Ker c_0
\subset C_{01}$. Therefore
$$
\Ker p \cong \Hom(e_*(\I),e_*(\I(M))),
$$
and by Yoneda, this is identified with $e_*(\I(M))(R) = M$. To
finish the proof, it remains to check that the map \eqref{p.rr} is
surjective --- that is, every element $r \in R$ lifts to some
regular endomorphism $r' \in \End(C_{01})$. Indeed, fix an element
$r$. Since the elementary extension $\phi = e_*(C_\idot)$ comes
from an elementary extension in $\B(R)$, we have
$$
r \circ \phi \cong \phi \circ r,
$$
where the compositions are given by the composition functors
\eqref{comp.l.ex}, \eqref{comp.r.ex}. Then a lifting $r'$ exists if
and only if we have
\begin{equation}\label{r.spl}
r \circ C_{01} \cong C_{01} \circ r \in \spl(r \circ \phi) \cong
\spl(\phi \circ r),
\end{equation}
where the compositions on the left-haand side are given by the
composition functors \eqref{comp.l.spl}, \eqref{comp.r.spl}. But the
extension $r \circ \phi$ is also of the form \eqref{phi.p.m.seq},
so that the existence of an isomorphism \eqref{r.spl} immediately
follows from the uniqueness statement of
Proposition~\ref{Z.spl.prop}.
\endproof

\subsection{The inversion theorem.}

By virtue of Lemma~\ref{regu.le} and Proposition~\ref{Z.spl.prop},
for any associative unital ring $R$ and $R$-bimodule $M$, we have a
well-defined map
\begin{equation}\label{el.sq}
\El(R,M) \to \Sq(R,M)
\end{equation}
that sends an elementary extension $\alpha$ to the corresponding
square-zero extension $R'$.

\begin{theorem}\label{sq.el.thm}
The maps \eqref{sq.el} and \eqref{el.sq} are inverse to each other.
\end{theorem}

For the proof of Theorem~\ref{sq.el.thm}, it is convenient to
introduce an additional piece of notation. For any square-zero
extension $R' \in \Sq(R,M)$, and any $R'$-module $V'$, the embedding
$M \subset R'$ induces a natural multiplication map
\begin{equation}\label{m.v.0}
M \otimes_{R'} V' \to R' \otimes_{R'} V' \cong V',
\end{equation}
and by \eqref{sq.eq}, $q^*V' = R \otimes_{R'} V'$ is the cokernel of
this map. Moreover, by adjunction, we have an identification $M
\otimes_{R'} V' \cong M \otimes_R q^*V'$. We will denote by
$m_{V'}$ the canonical map
\begin{equation}\label{m.v}
m_{V'}:M \otimes_R q^*V' \to V'
\end{equation}
induced by \eqref{m.v.0} via this identification. The module $V'$ is
of the form $q_*V$ for some $V \in R\amod$ if and only if
$m_{V'}=0$.

\proof[Proof of Theorem~\ref{sq.el.thm}.] Assume first given a
square-zero extension \eqref{sq.eq}, and consider the corresponding
elementary extension $C_\idot = \tau_{[0,1]}q^*P'_\idot$ in $\B(R)$,
where $P'_\idot$ is a projective resolution of the object $q_*(\I)
\in \B(R,R')$. Denote $C'_\idot = \tau_{[0,1]}P'_\idot$. Then by
adjunction, we have a natural map
$$
C'_\idot \to q_*C_\idot
$$
in the category $\B(R,R')$, and its composition with the projection
$C_\idot \to q_*(\I)$ is a quasiisomorphism. Then taking
$\wt{C}^l_\idot = C'_\idot$ in \eqref{DG.el.spl} provides a
splitting $C_{01}$ of the elementary extension $q_*C_\idot$. Denote
by $e':\Z \to R'$ the tautological map, and consider the induced
splitting $e'_*(C_{01})$ of the elementary extension
$e'_*(q_*(C_\idot)) \cong e_*(C_\idot)$. By uniqueness, it must be
isomorphic to the splitting provided by
Proposition~\ref{Z.spl.prop}. But by construction, the ring $R'$
acts on $e'_*(C_{01})$ by regular endomorphisms. Therefore it maps
to the square-zero extension provided by Lemma~\ref{regu.le}. Since
all maps between square-zero extensions are isomorphisms, this
proves that \eqref{el.sq} sends $C_\idot$ to $R'$.

Conversely, assume given an elementary extension $\langle
C_\idot,a,b \rangle$ with the corresponding splitting $\langle
C_{01},c_0,c_1 \rangle$ of Proposition~\ref{Z.spl.prop}, and let $R'
\subset \End(C_{01})$ be the square-zero extension provided by
Lemma~\ref{regu.le}. Moreover, consider the corresponding left DG
splitting $C^l_\idot$ of \eqref{dg.lr}. Then by definition, $R'$
acts on $C_{01}$ and $C^l_\idot$, so that $C_{01}$ actually defines
a splitting of the extension $q_*(C_\idot)$ in the category
$\B(R,R')$, and $C^l_\idot$ is a complex in $\B(R,R')$
quasiisomorphic to $q_*(\I)$. By the usual property of projective
resolutions, for any projective resolution $P'_\idot$ of $q_*(\I)$,
we then have a quasiisomorphism $P'_\idot \to C^l_\idot$, and it
induces a map
$$
\overline{P}_\idot = \tau_{[0,1]}q^*P'_\idot \to
\tau_{[0,1]}q^*C^l_\idot = q^*C^l_\idot.
$$
To finish the proof, it remains to show that $q^*C^l_\idot \cong
C_\idot$ --- indeed, then $\overline{P}_\idot$ and $C_\idot$ lie in
the same connected component of the category $\el(R,M)$, thus
correspond to the same element in $\El(R,M)$. In homological degree
$1$, we have $C^l_1 = q_*C_1$, so that $q^*C_1^l = q^*q_*C_1 \cong
C_1$. In degree $0$, $q^*C^l_0=q^*C_{01}$ is the cokernel of the map
$m_{C_{01}}$ of \eqref{m.v}. However, by construction, this maps
factors as
\begin{equation}\label{mult.fact}
\begin{CD}
M \otimes_R q^*C_{01} @>{\id \otimes (a \circ q^*c_0)}>> M \otimes_R
\I \cong \I(M) @>{c_1 \circ b}>> C_{01},
\end{CD}
\end{equation}
where the map $\id \otimes (a \circ q^*c_0)$ is surjective, and the
map $c_1 \circ b$ is injective. Therefore indeed $q^*C_0^l \cong
\Coker m_{C_{01}} \cong \Coker (c_1 \circ b) \cong C_0$.
\endproof

\section{Splittings and liftings.}\label{spl.sec}

\subsection{Modules.}\label{lft.mod.subs}

Assume given a unital associative ring $R$, an $R$-bimodule $M$, and
an elementary extension $C_\idot$ of $R$ by $M$, and let $R' \in
\Sq(R,M)$ be the corresponding square-zero extension. As we have
mentioned in the proof of Theorem~\ref{sq.el.thm}, the canonical
splitting $C_{01}$ of the elementary extension $e_*(C_\idot)$ is
naturally a left module over $R'$, so that it actually defines a
splitting of the elementary extension $q_*(C_\idot)$ in the category
$\B(R,R')$. We denote this splitting by the same letter $C_{01}$
since it is effectively the same object. We now want to explain how
$C_{01}$ helps to describe modules over $R'$ in terms of modules
over $R$.

We start with the following observation. By definition, for any two
associative unital rings $R$, $R'$, a functor $F \in \B(R,R')$ is
defined on the category $R\proj$ of finitely generated projective
left $R$-modules. However, we can extend its domain of definition in
the standard way. Namely, for any $R$-module $V$, the category
$I(V)$ of finitely presented $R$-modules $V_i$ equipped with a map
$i:V_i \to V$ is small and filtering, and we have
\begin{equation}\label{fl.lim}
V = \lim_{\overset{I(V)}{\to}}V_i.
\end{equation}
If $V$ is flat, then the full subcategory $J(V) \subset I(V)$
spanned by $\langle V_i,i \rangle$ with $V_i \in R\proj$ is cofinal
(see e.g.\ \cite[Lemma 3.5]{per} that assumes commutativiy of $R$
but does not use it). Therefore $J(V)$ is also filtering, and we can
replace the colimit in \eqref{fl.lim} with the colimit of the same
functor over $J(V)$. Therefore for any $F \in \B(R,R')$ and flat
$R$-module $V$, we can set
$$
F(V) = \lim_{\overset{J(V)}{\to}}F(V_i) \in
R'\amod,
$$
and if we denote by $R\flmod$ the category of flat $R$-modules, then
this is a well-defined functor from $R\flmod$ to
$R'\amod$. Moreover, the extension operation is functorial in $F$,
and since the categories of modules satisfy Grothendieck's axiom
$AB5$, the operation is exact.

In particular, in our setting, the elementary extension $C_\idot$
defines an elementary extension
\begin{equation}\label{c.v.el}
\begin{CD}
0 @>>> M \otimes_R V[1] @>>> C_\idot(V) @>>> V @>>> 0
\end{CD}
\end{equation}
of $R$-modules for any $V \in R\flmod$, and the splitting $C_{01}$
provides a splitting $C_{01}(V)$ of the corresponding elementary
extension $q_*(C_\idot(V))$ in $R'\amod$. Both are functorial in
$V$.

\begin{defn}
A {\em lifting} of a flat $R$-module $V \in R\flmod$ to a square-zero
extension $R'$ is a flat $R'$-module $V' \in R'\flmod$ equipped with
an isomorphism $q^*V' = R \otimes_{R'} V' \cong V$.
\end{defn}

All liftings of a given $V \in R\flmod$ form a groupoid in an obvious
way. We denote this groupoid by $\lift(V,R')$, and we denote by
$\Lift(V,R')$ the set of isomorphism classes of its objects.

\begin{prop}\label{lft.spl.prop}
For any $V \in R\flmod$, we have a natural equivalence of categories
$$
\lift(V,R') \cong \spl(C_\idot(V))
$$
between liftings of $V$ to $R'$ and splittings of the elementary
extension \eqref{c.v.el}.
\end{prop}

\proof{} For any lifting $V' \in \lift(R,V)$, the canonical map
\eqref{m.v} is injective and fits into a short exact sequence
$$
\begin{CD}
0 @>>> q_*(M \otimes_R V) @>{m_{V'}}>> V' @>>> q_*V @>>> 0.
\end{CD}
$$
Therefore every lifting is in partiticular an extension of $q_*V$ by
$q_*(M \otimes_R V)$, so that we have a natural functor
\begin{equation}\label{lft.ex}
\lift(R,V) \to \ex(q_*V,q_*(M \otimes_R V)).
\end{equation}
Moreover, for any extension $V' \in \ex(q_*V,q_*(M \otimes_R V))$,
the canonical map $m_{V'}$ factors through $q_*(M \otimes_R V)
\subset V'$ by means of a map
$$
\overline{m}_{V'}:q_*(M \otimes_R V) \to q_*(M \otimes_R V),
$$
and $V'$ comes from a lifting if and only if $\overline{m}_{V'} =
\id$ is the identity map. Therefore in particular, \eqref{lft.ex} is
a fully faithful embedding.

On the other hand, the functor $q_*$ induces a fully faithful
embedding
$$
q_*:\spl(C_\idot(V)) \to \spl(q_*C_\idot(V)).
$$
To describe its image, note that since the maps $m_{q_*C_0(V)}$ and
$m_{q_*C_1(V)}$ both vanish, the map $m_{C'}$ for any splitting $C'
\in \spl(q_*C_\idot(V))$ factors as
$$
\begin{CD}
M \otimes_R q^*C' @>{\id \otimes (a \circ q^*c_0)}>> M \otimes_R V
@>{\overline{m}_{C'}}>> M \otimes_R V @>{c_1 \circ b}>> C'(V),
\end{CD}
$$
where as in \eqref{mult.fact}, $\id \otimes (a \circ q^*c_0)$ is
surjective, $c_1 \circ b$ is injective, and $\overline{m}_{C'}$ is
some map. Then $C'$ lies in $q_*(\spl(C_\idot(V))) \subset
\spl(q_*C_\idot(V))$ if and only if $\overline{m}_{C'}=0$.

It remains to observe that the canonical splitting $C_{01} \in
\spl(q_*C_\idot)$ provides a pair of functors
\begin{align}
&\spl(q_*C_\idot) \to \ex(q_*V,q_*(M \otimes_R V)), \qquad C' \mapsto
C_{01} - C',\label{there}\\
&\ex(q_*V,q_*(M \otimes_R V)) \to \spl(q_*C_\idot), \qquad E \mapsto C_{01}
- E,\label{back.again}
\end{align}
where in the right-hand side, we have the difference functors
\eqref{diff.fu} and \eqref{sum.fu}, and these functors are mutually
inverse by \eqref{sum.diff}. Moreover, whenever $C' \cong C_{01} -
E$, we have
$$
\overline{m}_{C_{01}} = \overline{m}_{C'} + \overline{m}_E.
$$
But since by \eqref{mult.fact}, we have $\overline{m}_{C_{01}} =
\id$, this means that the functors \eqref{there} and
\eqref{back.again} induce an equivalence between the subcategories
$\lift(V,R') \subset \ex(q_*V,q_*(M \otimes_R V))$ and
$\spl(C_\idot(V)) \subset \spl(q_*C_\idot(V))$.
\endproof

\subsection{DG sections.}

We next want to extend the correspondence between liftings and
splitting obtained in Proposition~\ref{lft.spl.prop} to complexes of
$R$-modules and elementary DG extensions of
Definition~\ref{DG.ext.R.def}. In order to do this, we first need an
appropriate version of Proposition~\ref{Z.spl.prop}.

Fix an associative unital ring $R$, and consider the category
$\B_\idot(R,\Z)$ of \eqref{b.d.rr}. The forgetful functor
$e_*:C_\idot^{pf}(R) \to C_\idot(\Z)$ is admissible in the sense of
Definition~\ref{adm.def}, hence gives an object $e_*(\I_\idot) =
\I_\idot(R) \in \B_\idot(R,\Z)$. Denote
$$
K_\idot = \Cyl(R) \in C^{pf}_\idot(R),
$$
or explicitly, $K_0=K_1=R$, $K_i=0$ otherwise, the differential
$d:K_1 \to K_0$ given by the identity map. Then for any $m$, we have
the tautological map $\delta:K_\idot[m] \to K_\idot[m+1]$ given by
the composition
\begin{equation}\label{delta.eq}
\begin{CD}
K_\idot[m] @>{\alpha}>> R[m+1] @>{\beta}>> K_\idot[m+1]
\end{CD}
\end{equation}
of the tautological maps \eqref{cone.eq}. We have $\delta^2=0$, so
that $\langle K_\idot[i],\delta \rangle$ is a complex in the
additive category $C_\idot(R)$.

\begin{defn}\label{dg.sec.def}
A {\em DG section} $s_\idot$ of a surjective map $f:F_\idot \to
e_*(\I_\idot)$ in the category $\B_\idot(R,\Z)$ is a collection of
elements $s_m \in F_m(K_\idot[m-1])$, $m \in \Z$ such that
\begin{equation}\label{dg.sec.eq}
F_\idot(\delta)(s_m) = ds_{m+1}, \qquad f(s_m)=1 \in R =
e_*(\I_\idot)(K_\idot[m-1])_m)
\end{equation}
for any integer $m$.
\end{defn}

For any integer $m$, the complex $K_\idot[m-1]$ represents the
forgetful functor sending a complex $V_\idot$ to the set $V_m$, so
that a DG section $s_\idot$ represents by Yoneda a collection of
functorial maps
\begin{equation}\label{dg.s.yo}
\wt{s}_m:V_m \to F_m(V_\idot), \qquad V_\idot \in C^{pf}_\idot(R)
\end{equation}
such that $s_m \circ d = d \circ s_{m-1}$. A DG section does not
always exists. For example, let $F_\idot = \Cyl(e_*(\I_\idot)[-1])$
be the functor that sends a complex $V_\idot \in C_\idot(R)$ to the
underlying complex of $\Cyl(V_\idot)[-1]$, and let $f = \kappa_0$ be
the tautological map of \eqref{cone.eq}. Then it is easy to check
that $f:F \to e_*(\I_\idot)$ {\em does not} admit a DG section in
the sense of Definition~\ref{dg.sec.def}. However, as the following
result shows, it is easy to exclude situations of this kind.

\begin{lemma}\label{dg.sec.le}
Assume given surjective map $f:F_\idot \to e_*(\I_\idot)$ in the
category $\B_\idot(R,\Z)$, and assume that there exists an element
$t \in F_0(R)$ such that $f(t)=1 \in R$ and $dt=0$. Then $f$ admits
a DG section in the sense of Definition~\ref{dg.sec.def}.
\end{lemma}

\proof{} First of all, observe that once we have an element $t_m \in
F_m(R[m])$ for some integer $m$ such that $dt_m=0$ and $f(t_m)=1$,
it is trivial to construct an element $s_m \in F_{m+1}(K_\idot[m])$
such that $ds_m = F_\idot(\alpha)(t_m)$. Indeed, the admissible
functor $F_\idot$ must send the acyclic complex $K_\idot[m]$ to an
acyclic complex, and $dF_\idot(\alpha)(t_m) =
F_\idot(\alpha)(dt_m)=0$. Moreover, $ds_m=F_\idot(\alpha)(t_m)$
implies $df(s_m) = \alpha(f(t_m))=1$, and since the relevant
differential in $K_\idot[m]$ is the identity map, this yields
$f(s_m)=1$. We can then let $t_{m+1} = F_\idot(\beta)(s_m)$, and
repeat the procedure. By induction, we can therefore assume given
elements $s_n \in F_n(K_\idot[n-1])$ satisfying \eqref{dg.sec.eq}
for all $n \geq m$, with $d(s_m)=F_\idot(\alpha)(t_m)$.

To finish the proof, we need to apply induction in the other
direction. Thus it suffices to construct $s_{m-1} \in
F_m(K_\idot[m-1])$ and $t_{m-1} \in F_{m-1}(R[m-1])$ such that
$F_\idot(\beta)(s_{m-1}) = t_m$, $ds_{m-1} =
F_\idot(\alpha)(t_{m-1})$, $dt_{m-1}=0$, and $f(t_{m-1})=1$,
$f(s_{m-1})=1$.

To do this, note that since $F_\idot$ is admissible, it must send
the exact sequence \eqref{cone.eq} for the complex $R[m-1]$ to a
quasiexact sequence. Then as we have observed right after the
definition of a quasiexact sequence, we can always choose $s_{m-1}
\in F_m(K_\idot[m-1])$ such that $F_\idot(\alpha)(s_{m-1})=t_m$, and
$ds_{m-1} = F_\idot(\beta)(t_{m-1})$ for some $t_{m-1} \in
F_{m-1}(R[m-1])$. Moreover, $F_\idot(\beta)$ is injective. Therefore
$F_\idot(\beta)(dt_{m-1}) = d(ds_{m-1}) = 0$ implies
$dt_{m-1}=0$. It remains to notice that the maps $\alpha$ and
$\beta$ are isomorphisms in relevant degrees, so that
$\alpha(f(s_{m-1})) = f(t_m)=1$ implies $f(s_{m-1})=1$, and this in
turns implies $f(t_{m-1})=1$.
\endproof

\begin{corr}\label{ext.1.dg.corr}
Assume given a surjective map $f:F_\idot \to e_*(\I_\idot)$ in the
category $\B(R,\Z)$ satisfying the assumptions of
Lemma~\ref{dg.sec.le}, and assume further that
\begin{equation}\label{bnd.0}
F_\idot(C_{\leq m}^{pf}(R)) \subset C_{\leq m}(\Z)
\end{equation}
for any integer $m$. Then $f$ admits a one-side inverse
$g:e_*(\I_\idot) \to F_\idot$, $f \circ g = \id$.
\end{corr}

\proof{} By Lemma~\ref{dg.sec.le}, $f$ admits a DG section
$s_\idot$, and it suffices to prove that the corresponding maps
\eqref{dg.s.yo} are additive. As in the proof of
Proposition~\ref{Z.spl.prop}, it suffices to consider the universal
situation: we let $V_\idot = K_\idot[m-1] \oplus K_\idot[m-1]$ for
some integer $m$, with $v_1,v_2 \in R \oplus R = V_m$ given by $v_1 = 1
\oplus 0$, $v_2 = 0 \oplus 1$, and must prove that
$$
p_m = \wt{s}_m(v_1 + v_2) - \wt{s}_m(v_1) - \wt{s}_m(v_2)
$$
is equal to $0$. By the definition of a DG section, we have
$F_\idot(\delta \oplus \delta)(p_m) = d(p_{m+1})$, so that $q_m =
F_\idot(\alpha \oplus \alpha)(p_m) \in F_m(R[m],R[m])$ is
closed. But by definition, $q_m$ lies in the cross-effects component
$F_\idot(R[m],R[m])$. Since $F_\idot$ is admissible, this is an
acyclic complex, and by assumption, it is trivial in homological
degrees $> m$. Therefore $q_m=0$ for any $m$. But then $p_m$ is
closed for any $m$, and by the same argument, it also vanishes.
\endproof

\begin{prop}\label{dg.Z.spl.prop}
Assume given an $R^o$-module $W$ and a DG elementary extension
\begin{equation}\label{dg.2.spl}
\begin{CD}
0 @>>> \I_\idot(W) @>{b}>> C_\idot @>{a}>> e_*(\I_\idot) @>>> 0
\end{CD}
\end{equation}
in the category $\B_\idot(R,\Z)$ such that $C_\idot$ satisfies
\eqref{bnd.1} with $C_{\leq m+1}(R)$ replaced by $C_{\leq
  m+1}(\Z)$. Then the extension \eqref{dg.2.spl} admits a strict
right DG splitting, and such a strict right DG splitting is unique
up to an isomorphism.
\end{prop}

\proof{} As in the proof of Proposition~\ref{Z.spl.prop}, uniqueness
immediately follows from Corollary~\ref{ext.1.dg.corr}. Indeed, by
virtue of \eqref{bnd.1}, the difference $E_\idot = C^r_\idot \dst
\wt{C}^r_\idot$ between any two strict right DG splitting of
\eqref{dg.2.spl} satisfies the assumtions of
Corollary~\ref{ext.1.dg.corr}, so that we have a map $\I_\idot(W)
\oplus e_*(\I_\idot) \to E_\idot$ from the trivial extension
$\I_\idot(W) \oplus e_*(\I_\idot) \in \ex(e_*(\I_\idot),\I_\idot(W))
\subset \ex_\idot(e_*(\I_\idot),\I_\idot(W))$. Together with the
natural maps \eqref{df.ddf} and \eqref{dg.sum.diff}, it induces a
map $C^r_\idot \to \wt{C}^r_\idot$, and since all maps between
strict DG splittings are isomorphisms, $C^r_\idot$ is isomorphic to
$\wt{C}^r_\idot$.

To prove existence, note that by the definition of a DG elementaty
extension, the map $a:C_\idot(R) \to e_*(\I_\idot)(R)=R$ is an
isomorphism on homology in homological degree $0$. Therefore $a$
satisfies assumptions of Lemma~\ref{dg.sec.le} and admits a DG
section $s_\idot$. As in the proof of Proposition~\ref{Z.spl.prop},
fix such a section $s_\idot$ and consider the corresponding
functorial maps \eqref{dg.s.yo}. These are not necessarily additive,
but they do commute with the differentials. Now let
\begin{equation}\label{c.r.dg}
C^r_i(V_\idot) = C_i(V_\idot) \times V_{i-1}, \qquad i \in \Z,
\end{equation}
with the differential given by $d(c \times v)=dc + \wt{s}(v)$, and
with the maps $r$, $a^r$ given by $r(c) = c \times 0$, $a^r(c \times
v) = a(c)+v$. As in the proof of Proposition~\ref{Z.spl.prop}, in
order to turn $C^r_\idot$ into a strict right DG splitting, we need
to equip $C^r_i(V_\idot)$, $i \in \Z$ with functorial structures of
abelian groups compatible with the map $r$ and the differential, and
this amounts to constructing functorial maps
$$
c_i(-,-):V_i \times V_i \to C_{i+1}(V_\idot), \qquad i \in \Z
$$
satisfying \eqref{coc.eq.1}, \eqref{coc.eq.2}, \eqref{coc.eq.3}, and
\begin{equation}\label{coc.eq.4.dg.1}
dc_i(v_1,v_2) =
c_{i-1}(dv_1,dv_2)+\wt{s}_i(v_1+v_2)-\wt{s}_i(v_1)-\wt{s}_i(v_2)
\end{equation}
for any $v_1,v_2 \in V_i$, $i \in \Z$. As in the proof of
Proposition~\ref{Z.spl.prop}, to construct $c_i(-,-)$, it suffices
to consider the universal case $V_\idot = K_\idot[i-1] \oplus
K_\idot[i-1]$, with $v_1,v_2 \in R \oplus R = V_i$ given by $v_1 = 1
\oplus 0$, $v_2 = 0 \oplus 1$. Then we have
\begin{equation}\label{d.del}
c_{i-1}(dv_1,dv_2) = C_\idot(\delta \oplus
\delta)(c_{i-1}(v_1,v_2)),
\end{equation}
where $\delta$ is the canonical map \eqref{delta.eq}. Moreover, as
in \eqref{tw.eq}, the DG elementary extension \eqref{dg.2.spl}
generates an elementary extension $\rho_i(C_\idot) =
\rho(\tau^i(C_\idot))$ for any integer $i$, and the DG section
$s_\idot$ gives an element $s = C_\idot(\alpha)(s_i) \in C_i(R[i])$
such that $a(s)=1$. Then Proposition~\ref{Z.spl.prop} provides
functorial maps
$$
\wt{c}_i(-,-):V \times V \to C_{i+1}(V[i])
$$
satisfying \eqref{coc.eq.1}--\eqref{coc.eq.4} with $\wt{s}$ being
the map represented by $s$. Let us look for the elements
$c_i(v_1,v_2) \in C_i(K_\idot[i-1] \oplus K_\idot[i-1])$ with the
additional assumption
\begin{equation}\label{coc.eq.5}
C_\idot(\alpha \oplus \alpha)(c_i(v_1,v_2)) = \wt{c}_i(v_1,v_2).
\end{equation}
Then by virtue of \eqref{d.del} and \eqref{delta.eq},
\eqref{coc.eq.4.dg.1} can be rewritten as
\begin{equation}\label{coc.eq.4.dg.2}
dc_i(v_1,v_2) = C_\idot(\beta \oplus
\beta)(\wt{c}_{i-1}(v_1,v_2))+\wt{s}_i(v_1+v_2)-\wt{s}_i(v_1)-\wt{s}_i(v_2),
\end{equation}
and since $\wt{s}_i$ commutes with the differentials while the maps
$\wt{c}_\idot(-,-)$ satisfy \eqref{coc.eq.4}, the right-hand side is
annihilated by $d$. Again as in the proof of
Proposition~\ref{Z.spl.prop}, the condition \eqref{coc.eq.1} then
means that $c_i(-,-)$ takes values in the cross-effects component
$C_{i+1}(K_\idot[i-1],K_\idot[i-1])$, and \eqref{coc.eq.4.dg.2} with
\eqref{coc.eq.5} prescribe its differential and image
$C_\idot(\alpha \oplus \alpha)(c_i(-,-))$ in
$C_{i+1}(R[i],R[i])$. However, the map
$$
\begin{CD}
C_\idot(K_\idot[i-1],K_\idot[i-1]) @>{C_\idot(\alpha \oplus \alpha)}>>
C_\idot(R[i],R[i])
\end{CD}
$$
is a surjective maps of acyclic complexes, so that its kernel
$\overline{C}_\idot$ is acyclic, and moreover, by \eqref{bnd.1}, we
have $C_{i+2}(K_\idot[i-1])=0$, so that the differential
$d:\overline{C}_{i+1} \to \overline{C}_i$ in the acyclic complex
$\overline{C}_\idot$ is injective. Therefore \eqref{coc.eq.1},
\eqref{coc.eq.4.dg.2} and \eqref{coc.eq.5} uniquely determine
$c_\idot(-,-)$. As in the proof of Proposition~\ref{Z.spl.prop},
\eqref{coc.eq.2} and \eqref{coc.eq.3} now immediately follow from
this uniqueness.
\endproof

\subsection{Complexes.}\label{lft.comp.subs}

Assume now given an $R$-bimodule $M$, and a DG elementry extension
$\phi = \langle C_\idot,a,b \rangle \in \el_\idot(R,M)$. Then by
Theorem~\ref{sq.el.thm}, its restriction $r(\phi) \in \el(R,M)$ in
the sense of \eqref{e.r.ext} corresponds to a square-zero extension
$R'$ of $R$ by $M$. On the other hand, we have the forgetful functor
$e_*:\B_\idot(R) \to \B_\idot(R,\Z)$, and
Proposition~\ref{dg.Z.spl.prop} provides a canonical strict right DG
splitting $C^r_\idot$ of the extension $e_*(\alpha)$

\begin{lemma}\label{dg.regu.le}
With the notation above, the algebra $R'$ acts naturally on the
complex $C^r_\idot$, thus turing it into a strict right DG splitting
of the elementary extension $q_*(\alpha)$ in the category
$\B_\idot(R,R')$.
\end{lemma}

\proof{} As in the proof of Lemma~\ref{regu.le}, define regular
endomorphisms of the splitting $C^r_\idot$ as those that fit into a
diagram \eqref{r.r.r}, and denote by $R''$ the algebra of regular
automorphisms. Then the uniqueness statement of
Proposition~\ref{dg.Z.spl.prop} immediately implies that $R''$ is a
square-zero extension of $R$ by $M$.  Since restriction is
functorial, we have a natural map $R'' \to R'$, and since all maps
in $\sq(R,M)$ are isomorphisms, we have $R' \cong R''$.
\endproof

Now as in Subsection~\ref{lft.mod.subs}, we observe that for any
complex $V_\idot$ of flat $R$-modules, the category $J(V_\idot)$ of
complexes $V' \in C^{pf}_\idot(R)$ equipped with a map $i:V'_\idot
\to V_\idot$ is small and filtering, and we have
$$
V_\idot = \lim_{\overset{V'_\idot \in J(V_\idot)}{\to}} V'_\idot.
$$
Therefore setting
$$
C_\idot(V_\idot) = \lim_{\overset{V'_\idot \in J(V_\idot)}{\to}}
C_\idot(V'_\idot)
$$
extends the functor $C_\idot$ to a functor
\begin{equation}\label{c.dot.fl}
C_\idot:C^{fl}_\idot(R) \to C_\idot(R),
\end{equation}
where $C^{fl}_\idot(R) \subset C_\idot(R)$ is the full subcategory
spanned by complexes of flat $R$-modules, and for any complex
$V_\idot \in C^{fl}_\idot(R)$, we have a natural elementary
extension $\phi(V_\idot)$ of $V_\idot$ by $M \otimes_R V_\idot$ in
the category $C_\idot(R)$. Moreover, the right DG splitting
$C^r_\idot$ of the extension $q_*C_\idot$ provided by
Lemma~\ref{dg.regu.le} extends to a splitting of the extension of
\eqref{c.dot.fl} in the category of admissible functors from
$C^{fl}_\idot(R)$ to $C_\idot(R')$.

\begin{defn}
A {\em DG lifting} of a complex $V_\idot \in C^{fl}_\idot(R)$ to the
square-zero extension $R'$ is a complex $V'_\idot \in C^\hdot(R')$
equipped with a map $V'_\idot \to q_*V_\idot$ such that the
multiplication map $m_{V'_\idot}$ factors as
$$
\begin{CD}
q_*M \otimes_{R'} V'_\idot @>>> q_*(M \otimes_R V_\idot)
@>{\overline{m}}>> V'_\idot,
\end{CD}
$$
and the map $\overline{m}$ fits into a quasiexact sequence
$$
\begin{CD}
0 @>>> q_*(M \otimes_R V_\idot) @>{\overline{m}}>> V'_\idot @>>>
q_*V_\idot @>>> 0.
\end{CD}
$$
The category of DG liftings of $V_\idot$ to $R'$ is denoted
$\lift_\idot(V_\idot,R')$, and the set of its connected components
is denoted $\Lift_\idot(V_\idot,R')$.
\end{defn}

\begin{prop}\label{DG.lft.spl.prop}
For any complex $V_\idot \in C^{fl}_\idot(R)$ of flat $R$-modules,
there is a natural bijection between the sets
$\Lift_\idot(V_\idot,R')$ and $\Spl^l_\idot(\phi(V_\idot))$.
\end{prop}

\proof{} The proof is essentially the same as
Proposition~\ref{lft.spl.prop}. We consider the right DG splitting
$C^r_\idot$ of the extension $q_*(\phi)$ provided by
Lemma~\ref{dg.regu.le}, we let $C^{lr}_\idot$ be the left DG
splitting associated to it by the functor \eqref{lr.eq}, and we
define functors
$$
\begin{aligned}
&\spl_\idot^l(\phi(V_\idot)) \to \ex_\idot(q_*(V_\idot),q_*(M
\otimes_R V_\idot)),\\
&\ex_\idot(q_*(V_\idot),q_*(M \otimes_R V_\idot)) \to
\spl_\idot^l(\phi(V_\idot))
\end{aligned}
$$
by
$$
C'_\idot \mapsto C^{lr}_\idot - C'_\idot, \qquad E_\idot \mapsto
C^{lr}_\idot - E_\idot,
$$
where we use the difference functors \eqref{dg.sum.fu},
\eqref{dg.diff.fu}. Then we observe that by virtue of the canonical
maps \eqref{dg.sum.diff}, these functors induced mutually inverse
maps between $\Spl^l_\idot(\phi(V_\idot))$ and
$\Lift_\idot(V_\idot,R')$.
\endproof

\begin{remark}
The reader might notice that contrary to what we said in the
Introduction, we do not recover complexes in $C_\idot(R')$ as
algebras over a monad in $C_\idot(R)$ extending the endofunctor
$C_\idot$. The reason for this is that the monad in question is
essentially freely generated by $C_\idot$. Our splittings extend to
algebras over this monad canonically, and technically, there is no
reason to actually work this out: it adds a lot of complexity but
gives the same end result. However, if one want to go beyound
square-zero extensions, the structure of a monad or a comonad
becomes essential.
\end{remark}

\section{Multiplication.}\label{mult.sec}

\subsection{Multiplicative extensions.}

Now assume that our ring $R$ is commutative, so that the category
$R\proj$ of projective finitely generated $R$-modules is a tensor
category, with the unit object $R$.

\begin{defn}\label{mult.def}
For any commutative ring $R'$, a {\em multiplicative structure} on a
pointed functor $F \in \B(R,R')$ is given by a map $e:R' \to F(R)$
and a collection of functorial maps
$$
m(V,V'):F(V) \otimes_{R'} F(V') \to F(V \otimes_R V'), \qquad V,V' \in
R\proj,
$$
subject to obvious associativity and unitality conditions.
\end{defn}

For example, the tautological functor $\I \in \B(R)$ carries an
obvious multiplicative structure. More generally, for any
associative $R$-algebra $A$, we can treat $A$ as a diagonal
$R$-bimodule --- that is, take the same right $R$-action as the left
$R$-action, --- and then $\I(A)$ has a natural multiplicative
structure induced by the multiplication in $A$.

Conversely, for any associative $R$-algebra $A$ and any functor $F
\in \B(R)$ equipped with a multiplicative structure, $F(A)$ is
naturally an $R$-algebra. In particular, $F(R)$ is an $R$-algebra.
Then one can show that if $F = \I(M)$ is additive, then $F =
\I(F(R))$ (so that in particular, $M$ must be a diagonal
$R$-bimodule). Since we will not need it, we do not give a proof.

More generally, given a complex $F_\idot$ in $\B(R,R')$, we say that
a multiplicative structure on $F_\idot$ is given a map $e:R' \to
F_0(R)$ and a collection of functorial maps
$$
m:F_\idot(V) \otimes_{R'} F_\idot(V') \to F_\idot(V \otimes_R V'),
\qquad V,V' \in R\proj,
$$
again subject to associativity and unitality conditions.

Equivalently, one can consider the tensor product of $R$-modules as
a functor
$$
\mu:R\proj \times R\proj \to R\proj,
$$
and one defines a tensor structure on $\Fun(R\proj,R'\amod)$ by
setting
\begin{equation}\label{b.ten}
F \motimes F' = \mu_!(F \boxtimes_{R'} F'),
\end{equation}
where $\mu_!$ is the left Kan extension with respect to the functor
$\mu$. By adjunction, the product of pointed functors is pointed, so
that the subcategory $\B(R,R') \subset \Fun(R\proj,R'\amod)$ also
acquires a tensor structure. Then equipping $F \in \B(R,R')$ with a
multiplicative structure is equivalent to turning it into an algebra
object in $\langle \B(R,R'),\motimes \rangle$, and giving a
multiplicative structure on a complex $F_\idot$ is equivalent to
turning it into a DG algebra in $\B(R,R')$.

\begin{defn}
An elementary extension $\phi=\langle C_\idot,a,b\rangle$ of a
commutative ring $R$ by an $R$-bimodule $M$ is {\em multiplicative}
if $C_\idot$ is equipped with a multiplicative structure such that
$a:C_\idot \to \I$ is a multiplicative map.
\end{defn}

\begin{lemma}\label{diag.le}
Assume given a multiplicative elementary extension $\phi$ of a
commutative ring $R$ by an $R$-bimodule $M$. Then the bimodule $M$
is diagonal.
\end{lemma}

\proof{} By definition, since the extension $\phi = \langle
C_\idot,a,b \rangle$ is multiplicative, $C_0$ is an algebra object
in $\langle \B(R),\motimes \rangle$, and $C_1$ is a module over
$C_0$. By restriction, $\I(M) \subset C_1$ is also a module, so that
we have functorial maps
\begin{equation}\label{m.c}
m(V,V'):C_0(V) \otimes_R (M \otimes_R V') \to M \otimes_R (V
\otimes_R V'), \quad V,V' \in R\proj.
\end{equation}
But for any $m \in M \otimes_R V$, $c \in C_1(V')$, we have
$$
m(V,V')(dc \otimes m) = d(m(V,V')(c \otimes m)) + m(V,V')(c \otimes
dm) = 0,
$$
so that the $C_0$-action on $\I(M)$ vanishes on $d(C_1) \subset C_0$
and factors through $\I = C_0/d(C_1)$. Then if we take $V'=R$, the
map \eqref{m.c} is induced by a map
\begin{equation}\label{ch.sw}
V \otimes_R M = V \otimes_R (M \otimes_R R) \to M \otimes_R (V
\otimes_R R) = M \otimes_R V.
\end{equation}
This map is functorial in $V$, thus $\End_R(V)$-equivariant. If
$V=R$, both its source and its target are identified with $M$ by
unitality, and the map itself is the identity map. However, $R =
\End_R(R)$ acts via the left action on the bimodule $M$ on the
left-hand side of \eqref{ch.sw}, and via the right action in the
right-hand side.
\endproof

\begin{prop}\label{sq.mult.el.prop}
Assume given a square-zero extension $R' \in \sq(R,M)$ of a
commutative ring $R$ by an $R$-bimodule $M$, and assume that $R'$ is
commutative (in particular, $M$ is diagonal). Then the elementary
extension $\phi \in \el(R,M)$ that corresponds to $R'$ under the
equivalence of Theorem~\ref{sq.el.thm} can be chosen to be
multiplicative.
\end{prop}

\proof{} Note that the functor $q_*:R\proj \to R'\amod$ has a
natural multiplicative structure, so that the object $q_*(\I) \in
\B(R,R')$ is an algebra with respect to the tensor product
\eqref{b.ten}. On the other hand, the adjoint functor $q^*:\B(R,R')
\to \B(R)$ is a tensor functor. Thus by the construction of the
correspondence \eqref{sq.el}, it suffices to check that one can
choose a projective resolution $P_\idot$ of the algebra object
$q_*(\I) \in \B(R,R')$ that is a DG algebra with respect to the
tensor product \eqref{b.ten}. To do this, observe that by
adjunction, the tensor product $P \motimes P'$ of two projective
objects $P,P' \in \B(R,R')$ is projective. Therefore for any
projective $P \in \B(R,R')$, the tensor algebra
$$
T^\hdot P = \bigoplus_{i \geq 0}P^{\motimes i}
$$
is a projective object in $\B(R,R')$. Now we can construct a
multiplicative resolution of the algebra $q_*(\I)$ by the standard
inductive procedure.
\endproof

\subsection{Multiplicative splittings.}

We now want to invert Proposition~\ref{sq.mult.el.prop} by showing
that the correspondence \eqref{el.sq} sends multiplicative
elementary extensions to commutative square-zero extensions. It is
convenient to generalize the context slightly.

\begin{defn}\label{ele.mult.def}
An elementary extension \eqref{ele.ex} in a tensor abelian category
$\E$ is {\em multiplicative} if $A$ is an associative unital algebra
in $\E$, $C_\idot$ is a DG algebra, and $a:C_\idot \to A$ is an
algebra map. A {\em multiplicative splitting} of a multiplicative
elementary extension \eqref{ele.ex} in $\E$ is its splitting $\langle
C_{01},c_0,c_1 \rangle$ equipped with a DG algebra structure on the
complex $C^l_\idot$ of \eqref{dg.lr} such that the map $l:C^l_\idot
\to C_\idot$ is a DG algebra map.
\end{defn}

We note that for any splitting $C_{01}$, a DG algebra structure on
the complex $C^l_\idot$ is completely determined on its degree-$0$
part, that is, an algebra structure on $C_{01}$. In order for it to
define a multiplicative splitting, the map $c_0:C_{01} \to C_0$ must
be an algebra map, the left and right $C_{01}$-actions on itself
must factor through $c_0$, so that $C_{01}$ is a $C_0$-bimodule, and
the map $c_1:C_1 \to C_{01}$ must be a $C_0$-bimodule map.

For any commutative ring $R$, the forgetful functor $e_*:\B(R) \to
\B(R,\Z)$ is multiplicative. Therefore for any multiplicative
elemenetary extension $\phi=\langle C_\idot,a,b \rangle \in
\el(R,M)$ of $R$ by a diagonal $R$-bimodule $M$, the elementary
extension $e_*(\phi)$ is a multiplicative elementary extension in
$\B(R,\Z)$.

\begin{lemma}\label{mult.spl.le}
For any multiplicative elementary extension $\phi \in \el(R,M)$, the
canonical splitting $\langle C_{01},c_0,c_1 \rangle$ of the
elementary extension $e_*(\phi)$ provided by
Proposition~\ref{Z.spl.prop} is multiplicative in a natural way.
\end{lemma}

\proof{} Since $\phi=\langle C_\idot,a,b \rangle$ is multiplicative,
a natural choice of the element $s \in C_0(R)$ with $a(s)=1$ is
provided by the multiplicative structure on $C_0$, namely, by the
unity map $R \to C_0$. The corresponding map \eqref{s.yo} is then a
multiplicative map. By \eqref{p.times}, we have
$$
C_{01}(V) = C_1(V) \times V, \qquad V \in R\proj,
$$
with maps $c_0$, $c_1$ given by \eqref{c.0.1}. To turn $C_{01}$ into
an algebra in $\B(R,\Z)$, we need to construct multiplication maps
$$
- \cdot -:C_{01}(V) \times C_{01}(V') \to C_{01}(V \otimes_R V'),
\qquad V,V' \in R\proj
$$
that are associative, unital, and bilinear in each argument. We
note that by multiplicativity of the extension $C_\idot$, we have
\begin{equation}\label{sym}
c \cdot \delta(c') - \delta(c) \cdot c' = \delta(c \cdot c') = 0,
\qquad c,c' \in C_1(V),
\end{equation}
where $\delta:C_1(V) \to C_0(V)$ is the differential, and we let
\begin{equation}\label{m.vv}
\begin{aligned}
(c \times v) \cdot (c' \times v') &= (c \cdot \delta(c') +
c \cdot \wt{s}(v') + \wt{s}(v) \cdot c') \times (v \cdot v')\\
& = (\delta(c) \cdot c' + c \cdot \wt{s}(v') + \wt{s}(v) \cdot c') \times
(v \cdot v')
\end{aligned}
\end{equation}
for any $v \in V$, $v' \in V'$, $f \in C_1(V)$, $f' \in
C_1(V')$. This product is obviously unital, with $0 \times 1 \in
C_0(R)$ being the unit. To check that it is bilinear in the first
argument, we need to check that for any $v_1,v_2 \in V$, $c_1,c_2
\in C_1(V)$, we have
\begin{multline*}
(c_1 + c_2 + c(v_1,v_2)) \cdot (\delta(c') + \wt{s}(v')) + \wt{s}(v_1 +
v_2) \cdot c' =\\
= (c_1 + c_2) \cdot (\delta(c') + \wt{s}(v')) +
(\wt{s}(v_1) + \wt{s}(v_2)) \cdot c' + c(v_1 \cdot v',v_2 \cdot v'),
\end{multline*}
where $c(-,-)$ is the cocycle \eqref{coc}. By \eqref{coc.eq.4}, this
can be rewritten as
$$
c(v_1,v_2) \cdot \wt{s}(v') + c(v_1,v_2) \cdot \delta(c') -
\delta(c(v_1,v_2)) \cdot c' = c(v_1 \cdot v',v_2 \cdot v'),
$$
and taking into account \eqref{sym}, we further reduce it to
$$
c(v_1,v_2) \cdot \wt{s}(v') = c(v_1 \cdot v',v_2 \cdot v').
$$
By the same cross-effects argument as in the proof of
Proposition~\ref{Z.spl.prop}, it suffices to check this equality
after applying the differential $\delta$, and then it immediately
follows from \eqref{coc.eq.4} and multiplicativity of the map
$\wt{s}$.

The proof that the multiplication is bilinear in the second argument
is exactly the same --- all one has to do is to use the second of
the equivalent expressions in the right-hand side of \eqref{m.vv}
instead of the first one.

Now, by bilinearity, to prove that the multiplication \eqref{m.vv}
is associative, it suffices to check associativity separately for
elements of the form $c \times 0$ and $0 \times v$, and out of the
eight possibilities that arise, the only non-trivial one is the
equality
\begin{equation}\label{ass}
((c \times 0) \cdot (c' \times 0)) \cdot (c'' \times 0) = (c \times
0) \cdot ((c' \times 0) \cdot (c'' \times 0))
\end{equation}
for any $V,V',V'' \in R\proj$, $c \in C_1(V)$, $c' \in C_1(V')$,
$c'' \in C_1(V'')$. Then by \eqref{m.vv} and \eqref{sym}, the
left-hand side of \eqref{ass} is given by
$$
((c \cdot \delta(c')) \cdot \delta(c'')) \times 0 = ((\delta(c)
\cdot c') \cdot \delta(c'')) \times 0,
$$
and this clearly coincides with the right-hand side.

To finish the proof, it remains to prove that $c_0:C_{01} \to C_0$
is an algebra map, and $c_1:C_1 \to C_{01}$ is a $C_{01}$-module
map. This is a straightforward check that we leave to the reader; the
only non-trivial observation is that we have
$$
\delta(c) \cdot \delta(c') = \delta(c \cdot \delta(c'))
$$
for any $V,V' \in R\proj$, $c \in C_1(V)$, $c' \in C_1(V')$.
\endproof

\subsection{The inverse construction.}

We can now prove a converse to Proposition~\ref{sq.mult.el.prop}. We
keep the notation of last subsection.

\begin{prop}\label{el.mult.sq.prop}
Assume given a multiplicative elementary extension $\phi = \langle
C_\idot,a,b \rangle \in \el(R,M)$ of a commutative ring $R$ by a an
$R$-bimodule $M$, and let $R' \in \sq(R,M)$ be the square-zero
extension of $R$ by $M$ corresponding to $\phi$ by
\eqref{el.sq}. Then $R'$ is commutative, and the multiplicative
splitting $C_{01}$ of the induced extension $e_*(\phi)$ provided by
Lemma~\ref{mult.spl.le} extends to a multiplicative splitting of the
extension $q_*(\phi)$.
\end{prop}

\proof{} Since by assumption, $C_0$ is a multiplicative functor,
$C_0(R)$ is a unital associative ring, and the augmentation map
$C_0(R) \to R$ admits a natural splitting $R \to C_0(R)$. Since
$C_{01}$ is multiplicative, $C_{01}(R)$ is also a unital associative
ring, and $c_0:C_{01}(R) \to C_0(R)$ is a ring map. Consider the
subring $R'' = c_0^{-1}(R) \subset C_{01}(R)$, where $R \subset
C_0(R)$ is the image of the splitting $R \to C_0(R)$. Then $R''$ is
a square-zero extension of $R$ by $M$. Moreover, since $M$ is a
diagonal $R$-bimodule by Lemma~\ref{diag.le}, \eqref{m.vv}
immediately shows that $R''$ is commutative. Now, every $c \in
C_{01}(R)$ induces by \eqref{m.vv} a functorial map
\begin{equation}\label{times.c}
- \cdot c:C_{01}(V) \to C_{01}(V), \qquad V \in R\proj,
\end{equation}
and this map is obviously regular in the sense of
\eqref{r.r.r}. Therefore we have a natural ring map $C_{01}(R) \to
R'$. The restriction of this map to the square-zero extension $R''
\subset C_{01}(R)$ must be an isomorphism, so that $R' \cong R''$ is
indeed commutative. Moreover, the $R'$-action on $C_{01}$ can be
described by \eqref{times.c}, and this shows that the product in
$C_{01}$ is $R'$-linear, so that $C_{01}$ gives a multiplicative
splitting of $q_*(\phi)$.
\endproof

We can also prove a multiplicative refinement of
Proposition~\ref{lft.spl.prop}. Namely, assume given a flat
associative unital $R$-algebra $A$, and say that a {\em lifting}
$A'$ of $A$ to $R'$ is a flat associative unital $R'$-algebra $A'$
equipped with isomorphism $A' \otimes_{R'} R \cong A$. Denote the
groupoid of all liftings of $A$ to $R'$ by
$\lift^{\motimes}(A,R')$. On the other hand, note that since $A$ is
an algebra, $C_\idot(A)$ is a DG algebra in $R\amod$, and moreover,
it is a multiplicative elementary extension of $A$ by $M \otimes_R
A$ in the sense of Definition~\ref{ele.mult.def}. Denote by
$\spl^{\motimes}(C_\idot(A))$ the groupoid of multiplicative
splittings of this extension.

\begin{prop}\label{lft.mult.prop}
The equivalence of Proposition~\ref{lft.spl.prop} induces an
equivalence of categories
$$
\lift^{\motimes}(A,R') \cong \spl^{\motimes}(C_\idot(A)).
$$
\end{prop}

\proof{} By Proposition~\ref{el.mult.sq.prop}, the splitting
$C_{01}(A)$ of the extension $q_*C_\idot(A)$ is multiplicative. It
remains to observe that the difference functors \eqref{diff.fu} and
\eqref{sum.fu} are obviously compatible with multiplicative
structures.
\endproof

We note that the DG algebra $C_\idot(A)$ is a square-zero extension
of $A$ by the complex $(M \otimes_R A)[1]$, and multiplicative
splittings considered in Proposition~\ref{lft.mult.prop} are by
definition splittings of this extension in the category of DG
algebras over $R$. In particular, it is easy to construct an
obstruction to the existence of such a splitting: it is a certain
class in the third Hochschild cohomology group $HH^3(A,M \otimes_R
A)$. If the obstruction vanishes, then the set of isomorphism
classes of splittings is a torsor over $HH^2(A,M \otimes_R
A)$. Proposition~\ref{lft.mult.prop} then shows that the obstruction
theory for liftings of $A$ to $R'$ is exactly the same: there is an
obstruction lying in $HH^3(A,M \otimes_R A)$, and if it vanishes,
liftings are parametrized by a torsor over $HH^2(A,M \otimes_R
A)$. As far as we know, this seems to be a new result.

\subsection{DG extensions.}\label{dg.mult.subs}

A natural next thing to do would be to generalize the notion of a
multiplicative structure to DG elementary extensions of
Definition~\ref{DG.ext.R.def}, with the eventual goal of extending
Proposition~\ref{lft.mult.prop} to complexes of flat $R$-modules (or
equivalently, finding a multiplicative version of
Proposition~\ref{DG.lft.spl.prop}). By itself, a mutpliplicative
version of Definition~\ref{mult.def} presents no problems: all we
need to do is to add homological indices.

\begin{defn}\label{DG.mult.def}
For any two commutative rings $R$, $R'$, a {\em multiplicative
  structure} on an admissible functor $F_\idot \in \B_\idot(R,R')$
is given by a map $e:R' \to F_\idot(R)$ and a collection of
functorial maps
$$
m(V,V'):F_\idot(V_\idot) \otimes_{R'} F_\idot(V'_\idot) \to
F_\idot(V_\idot \otimes_R V'_\idot), \qquad V,V' \in C^{pf}_\idot(R),
$$
subject to obvious associativity and unitality conditions.
\end{defn}

Analogously, one can define multiplicative DG elementary extensions
and their multiplicative splitting. Moreover, for any multiplicative
DG extension $\phi_\idot \in \el_\idot(R,M)$, the restriction
$\rho(\phi_\idot) \in \el(R,M)$ is also multiplicative, so that the
associated square-zero extension $R' \in \Sq(R,M)$ must be a
commutative ring.

\medskip

But unfortunately, this is where the story ends, for the following
reason. The canonical strict DG splitting $C^r_\idot$ of the
extension $q_*(\phi_\idot)$ provided by Lemma~\ref{dg.regu.le} is a
right DG splitting, nor a left one; and the corresponding left DG
splitting $C^{rl}_\idot$ is not strict. On the other hand, if we
take the multiplicative splitting $C_{01}$ of
Proposition~\ref{el.mult.sq.prop}, then it is the corresponding
strict left DG splitting $C^l_\idot$ that is multiplicative. While
$C^l_\idot$ is quasiisomorphic to $C^{rl}_\idot$, it is not
isomorphic to it, and it seems that $C^{rl}_\idot$ does not have a
meaningful multiplicative structure.

Because of this, in the case of a general DG elementary extension,
we cannot proceed any further, and in particular, a DG version of
Proposition~\ref{lft.mult.prop} is beyond out reach.

\section{Cyclic powers.}\label{cycl.sec}

To finish the paper, we study in detail one particular DG elementary
extension $C_\idot$ of an arbitrary commutative ring $R$ annihilated
by a prime. We show how the general constructions work in this case,
and prove some additional results not available in the general
situation.

\subsection{Objects.}

Fix a commutative associative unital ring $R$ annhilated by a prime
$p > 1$. For any flat $R$-module $V \in R\flmod$, the $p$-th tensor
power $V^{\otimes_R p}$ is naturally a module over the group algebra
$R[\Z/p\Z]$ of the cyclic group $\Z/p\Z$, with the generator of the
group acting by the order-$p$ permutation $\sigma:V^{\otimes_R p}
\to V^{\otimes_R p}$. Moreover, the spaces of invariants and
coinvariants of the permutation $\sigma$ are related by a natural
trace map
\begin{equation}\label{c.v.tr}
\begin{CD}
\left(V^{\otimes_R p}\right)_\sigma @>{\tr}>> \left(V^{\otimes_R
  p}\right)^\sigma
\end{CD}
\end{equation}
given by $\tr = \id + \sigma + \dots + \sigma^{p-1}$.

Denote by $\wt{C}_\idot(V)$ the complex of $R$-modules obtained by
placing \eqref{c.v.tr} in homological degrees $0$ and $1$. Then the
homology of the complex $\wt{C}_\idot(V)$ has been computed in
\cite{ka3}. To state the answer, denote by $R^{(1)}$ the group $R$
considered as a module over itself via the Frobenius map --- that
is, let
$$
r \cdot r' = r^pr', \qquad r \in R, r' \in R^{(1)},
$$
and for any flat $R$-module $V$, let $V^{(1)} = R^{(1)} \otimes_R
V$. Note that the Frobenius map induces an $R$-module map $R \to
R^{(1)}$, hence a functorial map
\begin{equation}\label{fr}
\Fr:V \to V^{(1)}.
\end{equation}

\begin{lemma}[{{\cite[Lemma 6.9]{ka3}}}]\label{tr.le}
The complex $\wt{C}_\idot(V)$ fits into an exact sequence
$$
\begin{CD}
0 @>>> V^{(1)} @>{\psi}>> \left(V^{\otimes_R p}\right)_\sigma
@>{\tr}>> \left(V^{\otimes_R p}\right)^\sigma @>{\wh{\psi}}>>
V^{(1)} @>>> 0,
\end{CD}
$$
and this sequence is functorial with respect to $V$.\endproof
\end{lemma}

Explicitly, the map $\psi$ in Lemma~\ref{tr.le} is given by
$$
\psi(r \cdot v) = rv^{\otimes_R p}
$$
(this is additive modulo $\Im(\id - \sigma)$). If $V \in R\proj$,
then $\wh{\psi}$ is the dual map, and for a general flat $R$-module
$V$, it is obtained by taking the filtered colimit. In any case,
both $\psi$ and $\wh{\psi}$ are completely functorial.

\begin{defn}\label{c.def}
For any $V \in R\flmod$, the complex $C_\idot(V)$ is obtained by the
pullback square
\begin{equation}\label{sq.c.c}
\begin{CD}
C_\idot(V) @>{a}>> V\\
@VVV @V{\Fr}VV\\
\wt{C}_\idot(V) @>{\wt{a}}>> V^{(1)},
\end{CD}
\end{equation}
where $\wt{a}$ is the natural map induced by the map $\wh{\psi}$ of
Lemma~\ref{tr.le}, and $\Fr$ is the Frobenius map \eqref{fr}.
\end{defn}

By definition, $C_\idot$ is a complex of pointed functors from
$R\proj$ to $R\amod$, thus a complex in the category
$\B(R)$. Moreover, $C_\idot(V)$ comes equipped with a functorial map
$a:C_\idot(V) \to V$, and a functorial map $b:V^{(1)}[1] \to
C_\idot(V)$ is induced by the map $\psi$ of
Lemma~\ref{tr.le}. Altogether, $\langle C_\idot,a,b \rangle$ define
an elementary extension of $R$ by the diagonal bimodule $R^{(1)}$
that we denote by
\begin{equation}\label{phi.df}
\phi \in \el(R,R^{(1)}).
\end{equation}
Therefore by \eqref{el.sq}, there is a square-zero extension $R' \in
\el(R,R^{(1)})$ associated to the elementary extension $\phi$.

Moreover, for any flat $R$-modules $V_1,V_2 \in R\flmod$, we have
natural maps
$$
\begin{aligned}
&V_1^{\otimes_R p} \otimes_R \left(V_2^{\otimes_R p}\right)^\sigma \to
V_1^{\otimes_R p} \otimes_R V_2^{\otimes_R p} \cong (V_1 \otimes_R
V_2)^{\otimes_R p},\\
&\left(V_1^{\otimes_R p}\right)^\sigma \otimes_R V_2^{\otimes_R p} \to
V_1^{\otimes_R p} \otimes_R V_2^{\otimes_R p} \cong (V_1 \otimes_R
V_2)^{\otimes_R p},
\end{aligned}
$$
and these maps are $\sigma$-equivariant, thus descend to functorial
maps
$$
\begin{aligned}
&\left(V_2^{\otimes_R p}\right)^\sigma \otimes_R
  \left(V_2^{\otimes_R p}\right)^\sigma \to \left((V_1 \otimes_R
  V_2)^{\otimes_R p}\right)^\sigma,\\
&\left(V_2^{\otimes_R p}\right)_\sigma \otimes_R
  \left(V_2^{\otimes_R p}\right)^\sigma \to \left((V_1 \otimes_R
  V_2)^{\otimes_R p}\right)_\sigma,\\
&\left(V_2^{\otimes_R p}\right)^\sigma \otimes_R
  \left(V_2^{\otimes_R p}\right)_\sigma \to \left((V_1 \otimes_R
  V_2)^{\otimes_R p}\right)_\sigma
\end{aligned}
$$
compatible with the trace maps $\tr$. Thus the complex of functors
$\wt{C}_\idot(-)$ has a natural multiplicative structure. Then so
does the complex $C_\idot(-)$ of Definition~\ref{c.def}, so that the
elementary extension $\phi$ of \eqref{phi.df} is multiplicative. By
Proposition~\ref{el.mult.sq.prop}, the square-zero extension $R'$
must then be commutative.

As it happens, one commutative square-zero extension of $R$ by
$R^{(1)}$ is very well-known: this is the ring $W_2(R)$ of second
Witt vectors of $R$. Let us prove that this is exactly what we get
from the extension \eqref{phi.df}.

\begin{prop}
For any commutative ring $R$ annihilated by a prime $p$, the
square-zero extension $R' \in \Sq(R,R^{(1)})$ corresponding to the
elementary extension \eqref{phi.df} is isomorphic to the second Witt
vectors ring $W_2(R)$.
\end{prop}

\proof{} Note that by definition, we have $C_0(R) \cong R$, so that
by the same argument as in the proof of
Proposition~\ref{el.mult.sq.prop}, we have $R' \cong C_{01}(R)$,
where $C_{01}$ is the multiplicative splitting of the extension
$C_\idot$ provided by Lemma~\ref{mult.spl.le}. We also have $C_1(R)
\cong R^{(1)}$, so that as a set, we have
$$
R' \cong R \times R^{(1)}.
$$
The additive structure is given by Proposition~\ref{Z.spl.prop} ---
we have
\begin{equation}\label{sum.witt}
(x_0 \times x_1) + (y_0 \times y_1) = (x_0 + y_0) \times (x_1 +
y_1 + c(x_0,y_0))
\end{equation}
for any $x_0,y_0 \in R$, $x_1,y_1 \in R^{(1)}$, where $c(-,-)$ is
the cocycle \eqref{coc}. To compute $c(-,-)$ explicitly, note that
the multiplicative splitting map \eqref{s.yo} is given by $\wt{s}(v)
= v^{\otimes p}$, and consider the universal situation $V = Rv_1
\oplus Rv_2$. Then \eqref{coc.eq.4} reads as
$$
(\id + \sigma + \dots + \sigma^{p-1})c(v_1,v_2) = \sum_{s \in
  \overline{S}} \mu_s(v_1,v_2),
$$
where the sum is over the set $\overline{S}$ of all degree-$p$
non-commutative monomials in $v_1$ and $v_2$ except for
$v_1^{\otimes p}$ and $v_2^{\otimes p}$. In other words, we have
$\overline{S} = S^p \setminus S$, where $S = \{1,2\}$ is the set of
indices, and $S \subset S^p$ is the diagonal. The permutation action
of $G = \Z/p\Z$ on $\overline{S}$ is free, and to write down an
explicit formula for $c(v_1,v_2)$, it suffices to chose a splitting
$\kappa:\overline{S}/G \to \overline{S}$ of the quotient map
$\overline{S} \to \overline{S}/G$. We then have
\begin{equation}\label{c.exp}
c(v_1,v_2) = \sum_{s \in \overline{S}/G}\mu_{\kappa(s)}(v_1,v_2) \in
V^{\otimes_R p}_\sigma,
\end{equation}
and this does not depend on the choice of $\kappa$. Projecting back
onto $V=R$, and noting that $R$ is commutative, we see that
$$
c(x,y) = \sum_{1 \leq i \leq p-1}\frac{1}{p}\binom{p}{i}x^iy^{p-i},
\qquad x,y \in R,
$$
so that \eqref{sum.witt} becomes
\begin{equation}\label{sum.witt.bis}
(x_0 \times x_1) + (y_0 \times y_1) = (x_0 + y_0) \times \left(x_1 +
  y_1 + \sum_{1 \leq i \leq
    p-1}\frac{1}{p}\binom{p}{i}x_0^iy_0^{p-i}\right).
\end{equation}
As for the product in $R' = C_{01}(R)$, then by
Lemma~\ref{mult.spl.le}, it is given by \eqref{m.vv}, and in our
situation, it reads as
\begin{equation}\label{prod.witt}
(x_0 \times x_1) \cdot (y_0 \times y_1) = x_0y_0 \times
  (x_0^py_1+x_1y_0^p).
\end{equation}
It remains to recall that $R$ is assumed to be annihilated by $p$,
and to notice that modulo $p$, \eqref{sum.witt.bis} and
\eqref{prod.witt} are exactly the standard formulas for the sum and
product of Witt vectors.
\endproof

\subsection{Complexes.}\label{cycl.dg.subs}

The elementary extension $\phi \in \el(R,R^{(1)})$ defined by the
complex $C_\idot$ of Definition~\ref{c.def} also has a very natural
lifting to a DG elementary extension. To describe it, we first need
another interpretation of the complex $\wt{C}_\idot(V)$. Recall that
for any $R[\Z/p\Z]$-module $E$, the {\em Tate cohomology complex}
$\vC_\idot(\Z/p\Z,E)$ is obtained by taking $\vC_i(\Z/p\Z,E) = E$
for any integer $i$, with the differential $d_i:\vC_i(\Z/p\Z,E) \to
\vC_{i-1}(\Z/p\Z,E)$ given by
$$
d_i = \begin{cases} \id - \sigma, &\quad i=2j,\\ \tr, &\quad i=2j+1.
\end{cases}
$$
In other words, $\vC_\idot(\Z/p\Z,E)$ is the $2$-periodic complex
\begin{equation}\label{perio}
\begin{CD}
@>{\id - \sigma}>> E @>{\tr}>> E @>{\id-\sigma}>> E @>{\tr}>> E
@>{\id-\sigma}>>.
\end{CD}
\end{equation}
Then one immediately observes that we have a functorial
identification
\begin{equation}\label{c.tate}
\wt{C}_\idot(V) \cong \tau_{[0,1]}\vC_\idot(\Z/p\Z,V^{\otimes_R p})
\end{equation}
for any flat $R$-module $V$. To generalize \eqref{c.tate} to
complexes, one needs to find appropriate versions of the Tate
cohomology complex $C_\idot(\Z/p\Z,-)$ and of the truncation functor
$\tau_{[01]}$.

The former is easy: for any complex of $R[\Z/p\Z]$-modules
$E_\idot$, \eqref{perio} is naturally a bicomplex, and we denote by
$\vC_\idot(\Z/p\Z,E_\idot)$ its total complex.

For the latter, we use the filtered truncation functors introduced
in \cite[Section 1.3]{ka3}. Namely, for any abelian category $\E$,
denote by $CF_\idot(\E)$ the category of complexes in $\E$ equipped
with a descreasing filtration $F^\hdot$ numbered by all integers.

\begin{defn}\label{trunc.def}
For any integer $n$ and any abelian category $\E$, the filtered
truncation functors $\tau^F_{\leq n},\tau^F_{\geq n}:CF_\idot(\E)
\to CF_\idot(\E)$ are given by
$$
\begin{aligned}
\tau^F_{\geq n}E_i &= d^{-1}(F^{n+1-i}E_{i-1}) \cap F^{n-i}E_i
\subset E_i,\\
\tau^F_{\leq n}E_i &= E_i/(F^{n+1-i}E_i + d(F^{n-i}E_{i+1}))
\end{aligned}
$$
for any $E_\idot \in CF_\idot(\E)$.
\end{defn}

This is essentially \cite[(1.8)]{ka3}, with a difference of
notation: $\tau^F_{\geq n}$ is $\tau^n$ of \cite{ka3}, and
$\tau^F_{\leq n}E_\idot$ is the quotient
$E_\idot/\beta^nE_\idot$. For any filtered complex $E_\idot$ and any
$n$, we have natural maps $\tau^F_{\geq n}E_\idot \to E_\idot$,
$E_\idot \to \tau^F_{\leq n}$. Given two integers $n \leq m$, we
denote $\tau^F_{[n,m]}E_\idot = \tau^F_{\geq n}\tau^F_{\leq
  m}E_\idot \cong \tau^F_{\leq m}\tau^F_{\geq n}E_\idot$. If the
filtration $F^\hdot$ on a complex $E_\idot$ is termwise-split, then
we have natural identifications
\begin{equation}\label{gr.tau}
\begin{aligned}
\gr^i_F(\tau^F_{\geq n}E_\idot) &\cong \tau_{\geq n+i}\gr^i_FE_\idot,\\
\gr^i_F(\tau^F_{\leq n}E_\idot) &\cong \tau_{\leq n+i}\gr^i_FE_\idot,\\
\gr^i_F(\tau^F_{[n,m]}E_\idot) &\cong \tau_{[n+i,m+i]}\gr^i_FE_\idot
\end{aligned}
\end{equation}
for any integer $i$. In particular, $\tau^F_{\geq n}$, $\tau^F_{\leq
  n}$ and $\tau^F_{[n,m]}$ send filtered quasiisomorphisms to
filtered quasiisomorphism --- in effect, inverting filtered
quasiisomorphisms of filtered complexes gives the filtered derived
category $\DF(\E)$, and $\tau^F_{\leq n}$, $\tau^F_{\geq n}$ are
truncation functors with respect to a natural $t$-structure on
$\DF(\E)$.

We now note that any filtration $F^\hdot$ on a complex $E_\idot$ of
$R[\Z/p\Z]$-modules induces a filtration on the Tate complex
$\vC_\idot(\Z/p\Z,E_\idot)$. We denote by $\overline{F}^\hdot$ the
stupid filtration on $E_\idot$ --- we recall that by definition, it
is given by
$$
\overline{F}^iE_j = \begin{cases} E_j, &\quad i+j \geq 0,\\0, &\quad i+j < 0,
\end{cases}
\qquad i,j \in \Z.
$$
Then we define a filtration $F^\hdot$ on $E_\idot$ by rescaling
$\overline{F}^\hdot$ by $p$ --- that is, we let $F^iE_\idot =
\overline{F}^{ip}E_\idot$. With these definitions in mind, the
following generalization of Lemma~\ref{tr.le} has been obtained in
\cite{ka3}.

\begin{lemma}[{\cite{ka3}}]\label{dg.tr.le}
Assume given a complex $V_\idot$ of flat $R$-modules, and consider
the $p$-fold tensor power $E_\idot = V_\idot^{\otimes_R p}$. Then for any
integer $i$ not divisible by $p$, the Tate complex
$\vC_\idot(\Z/p\Z,E_i)$ is acyclic. Moreover, equip $E_\idot$ with
with the $p$-th rescaling $F^\hdot$ of the stupid filtration
$\overline{F}^\hdot$, and consider the induced filtration on the Tate
complex $\vC_\idot(\Z/p\Z,E_\idot)$. Then for any integer $i$, we
have an functorial isomorphism
\begin{equation}\label{dg.tr.eq}
\tau^F_{[i,i]}\vC_\idot(\Z/p\Z,E_\idot) \cong
V_\idot^{(1)}[i],
\end{equation}
this isomorphism is functorial in $V_\idot$, and the induced
filtration $F^\hdot$ on $V_\idot^{(1)}[i]$ is the stupid filtration
shifted by $i$.
\end{lemma}

\proof{} This is \cite[Proposition 6.10]{ka3} (the first and the
last claims are a part of the definition of a ``tight'' complex
given in \cite[Definiton 5.2]{ka3}).
\endproof

As an application of Lemma~\ref{dg.tr.le}, for any complex $V_\idot
\in C^{fl}_\idot(R)$ of flat $R$-modules, let us denote by
\begin{equation}\label{wt.dg.c}
\wt{C}_\idot(V_\idot) = \tau^F_{[0,1]}\vC_\idot(\Z/p\Z,V^{\otimes_R
  p}_\idot)
\end{equation}
the filtered truncation of the Tate cohomology complex
$\vC_\idot(\Z/p\Z,V^{\otimes_R p}_\idot)$. Then Lemma~\ref{dg.tr.le}
and \eqref{gr.tau} show that the complex $\wt{C}_\idot(V_\idot)$
fits into a natural functorial quasiexact sequence
$$
\begin{CD}
0 @>>> V^{(1)}_\idot[1] @>{b}>> \wt{C}_\idot(V_\idot) @>{\wt{a}}>>
V^{(1)}_\idot @>>>
0,
\end{CD}
$$
where as before, we let $V_\idot^{(1)} = R^{(1)} \otimes_R V_\idot$.
If we now define $C_\idot(V_\idot)$ by the pullback square
\eqref{sq.c.c}, then we have a functorial quasiexact sequence
\begin{equation}\label{dg.c.v}
\begin{CD}
0 @>>> V^{(1)}[1] @>{b}>> C_\idot(V_\idot) @>{a}>>
V_\idot @>>> 0.
\end{CD}
\end{equation}

\begin{lemma}\label{cycl.dg.le}
The sequence \eqref{dg.c.v} defines a DG elementary extension
$\phi_\idot$ of $R$ by $R^{(1)}$.
\end{lemma}

\proof{} It is immediately clear from \eqref{gr.tau} that the
functor $C_\idot = C_\idot(-)$ satisfies \eqref{bnd.1}, so it
remains to prove that it is admissible in the sense of
Definition~\ref{adm.def}. As we have remarked after stating the
condition \eqref{bnd.1}, the quasiexact sequence \eqref{dg.c.v}
insures that it suffices to prove that $C_\idot$ sends
termwise-split injections to termwise-split injections, and
termwise-split surjections to termwise-split surjections. In fact,
since the square \eqref{sq.c.c} is Cartesian, it suffices to prove
the same for the functor $\wt{C}_\idot$ of \eqref{wt.dg.c}. We will
do the injections --- the argument for surjections is exactly the
same.

Assume given a map $f:V_\idot \to V_\idot'$ of complexes of flat
$R$-modules, and denote
$$
E_\idot = V_\idot^{\otimes_R p}, \qquad E'_\idot = V_\idot^{'
  \otimes_R p}.
$$
Assume that $f$ is a termwise-split injection. Then the $p$-th
tensor power $f^{\otimes p}:E_\idot \to E'_\idot$ is a
termwise-split injection of complexes of $R[\Z/p\Z]$-modules, so
that for any integers $i$, $j$, the map
\begin{equation}\label{j.gr}
f^{\otimes p}:\tau_{[i,i+1]}\vC_\idot(\Z/p\Z,E_j) \to
\tau_{[i,i+1]}\vC_\idot(\Z/p\Z,E'_j)
\end{equation}
is an injection. Fix an integer $i$, and consider the associated
graded quotient $\gr^i_FE_\idot$. By the definition of the
filtration $F^\hdot$, this quotient has amplitude $[(i-1)p+1,ip]$,
and its associated graded quotients with respect to the stupid
filtration $\overline{F}^\hdot$ are the complexes $E_j[j]$,
$(i-1)p < j \leq ip$. Then the Tate complex
$\vC_\idot(\Z/p\Z,\gr^i_FE_\idot)$ is an iterated extension of
shifts of Tate complexes $\vC_\idot(\Z/p\Z,E_j)$, $(i-1)p < j \leq
ip$. But by Lemma~\ref{dg.tr.le}, all these complexes are acyclic,
except possibly for the one corresponding to $j=ip$. Then
applying inductively Lemma~\ref{triv.le}, we conclude that
$$
\gr^j_{\overline{F}}\tau_{[i,i+1]}\vC_\idot(\Z/p\Z,\gr^i_FE_\idot)
\cong \tau_{[i,i+1]}\vC_\idot(\Z/p\Z,E_j)[j], \quad (i-1)p < j \leq
ip.
$$
Moreover, we have the same identification for $E'_\idot$, and since
the maps \eqref{j.gr} are injections, we conclude that the
map
$$
\tau_{[i,i+1]}\vC_\idot(\Z/p\Z,\gr^i_FE_\idot) \to
\tau_{[i,i+1]}\vC_\idot(\Z/p\Z,\gr^i_FE'_\idot)
$$
induced by $f$ is an injection. Collecting these maps for all $i$
and applying \eqref{gr.tau}, we further conclude that the map
$$
\tau_{[0,1]}^F\vC_\idot(\Z/p\Z,E_\idot) \to
\tau_{[0,1]}^F\vC_\idot(\Z/p\Z,E'_\idot)
$$
induced by $f$ is an injection, and this is exactly what we had to prove.
\endproof

\subsection{Multiplication.}\label{mult.cycl.subs}

It is clear from \eqref{c.tate} and \eqref{gr.tau} that the
restriction $\rho(\phi_\idot) \in \el(R,R^{(1)})$ of the DG
elementary extension $\phi_\idot \in \el_\idot(R,R^{(1)})$ provided
by Lemma~\ref{cycl.dg.le} coincides with the elementary extension
$\phi$ of \eqref{phi.df}. As it happens, one can also extend the
multiplicative structure on $\phi$ to a multiplicative structure on
$\phi_\idot$. To do this, we need to recall a more invariant
definition of Tate cohomology.

For any finite group $G$ and any bounded complex $E_\idot$ of
$R[G]$-modules, the {\em Tate cohomology groups}
$\vH^\hdot(G,E_\idot)$ are given by
$$
\vH^\hdot(G,E) = \Ext^\hdot_{\D^b(R[G])/\D^{pf}(R[G])}(R,E_\idot),
$$
where $R$ is the trivial $R[G]$-module, and the $\Ext$-groups are
computed in the quotient of the bounded derived category
$\D^b(R[G])$ of $R[G]$-modules by its full subcategory
$\D^{pf}(R[G]) \subset \D^b(R[G])$ spanned by perfect complexes of
$R[G]$-modules (equivalently, $\D^{pf}(R[G]) \subset \D^b(R[G])$ is
spanned by compact objects in the triangulated category
$\D(R[G])$). In particular, $\vH^\hdot(G,R)$ is always an algebra.

In the case $G=\Z/p\Z$, the Tate cohomology complex
$\vC_\idot(\Z/p\Z,E)$ of \eqref{perio} computes exactly the groups
$\vH^\hdot(\Z/p\Z,E)$. Unfortunately, the multiplication in
$\vH^\hdot(\Z/p\Z,R)$ does not lift to a DG algebra structure on
$\vC_\idot(\Z/p\Z,R)$. However, this problem can be solved by
changing the complex.

Namely, consider a finite group $G$, choose a projection resolution
$P_\idot$ of the trivial $\Z[G]$-module $\Z$, and let
$\overline{P}_\idot$ be the cone of the augmentation map $P_\idot
\to \Z$. Then it is easy to show (see e.g.\ \cite[Subsection
  7.2]{ka-ma} but the claim is completely standard) that for bounded
complex $E_\idot$ of $R[G]$-modules, we have a natural
identification
\begin{equation}\label{dg.tate.df}
\vH^\hdot(G,E_\idot) \cong \lim_{\overset{i}{\to}}H^\hdot(G,E_\idot
\otimes F^{-i}\overline{P}_\idot),
\end{equation}
where $H^\hdot(G,-) = \Ext^\hdot_{\D(R[G])}(R,-)$ is the cohomology
of the group $G$, and $F^{-i}\overline{P}_\idot$ is the $(-i)$-th
term of the stupid filtration on the complex
$\overline{P}_\idot$. Moreover, the right-hand side of
\eqref{dg.tate.df} is canonically independent of the choice of a
resolution $P_\idot$ --- indeed, for any two resolutions $P_\idot$,
$P'\idot$, we have natural maps $\overline{P}_\idot \to
\overline{P}_\idot \otimes \overline{P}'_\idot$,
$\overline{P}'_\idot \to \overline{P}_\idot \otimes
\overline{P}'_\idot$, and both maps induce isomorphisms in the
right-hand side of \eqref{dg.tate.df}. If one wants to represent
Tate cohomology by an explicit functorial complex, one also needs to
represent the usual cohomology $H^\hdot(G,E_\idot)$ by such a
complex; the standard way to do is is to choose an injective
resolution $I^\hdot$ of the trivial $\Z[G]$-module $\Z$, and
consider the cohomology complex
\begin{equation}\label{c.do}
C^\hdot(G,E_\idot) = \left(E_\idot \otimes I^\hdot\right)^G.
\end{equation}
Altogether, the following sums up the situation.

\begin{defn}
\begin{enumerate}
\item {\em Resolution data} for a finite group $G$ is a pair $\nu =
  \langle P_\idot,I^\hdot \rangle$ of a projective resolution
  $P_\idot$ and an injective resolution $I^\hdot$ of the trivial
  $\Z[G]$-module $\Z$.
\item For any associative unital ring $R$, any bounded complex
  $E_\idot$ of $R[G]$-modules, and any resolution data $\nu$, the
  {\em Tate cohomology complex} of $G$ with coefficients in
  $E_\idot$ is given by
\begin{equation}\label{c.phi.tate}
\vC_\idot(G,\nu,E_\idot) = \left(E_\idot \otimes \overline{P}_\idot
\otimes I_\idot\right)^G,
\end{equation}
where $\overline{P}_\idot$ is the cone of the augmentation map $P_\idot
\to \Z$.
\end{enumerate}
\end{defn}

Then for any $G$ and $E_\idot$, and for any choice of the resolution
data $\nu$, we obviously have
$$
\vC_\idot(G,\nu,E_\idot) \cong \lim_{\overset{i}{\to}}C^\hdot(G,E_\idot
\otimes F^{-i}\overline{P}_\idot),
$$
where $C^\hdot(G,-)$ in the right-hand side is the complex
\eqref{c.do}, and by \eqref{dg.tate.df}, $\vC_\idot(G,\nu,E_\idot)$
computes the Tate cohomology $\vH_\idot(G,E_\idot)$ and does not
depend on $\nu$ up to a canonical quasiisomorphism.

Alternatively, given resolution data $\nu = \langle P_\idot,I^\hdot
\rangle$, one can consider the composition $P_\idot \to \Z \to
I^\hdot$ of the augmentation maps, and let $\wt{P}_\idot$ be its
cone. Then for any bounded complex $E_\idot$ of $R[G]$-modules, one
can consider the complex
\begin{equation}\label{c.dash.tate}
\vC'_\idot(G,\nu,E_\idot) = \left(E_\idot \otimes \wt{P}_\idot\right)^G.
\end{equation}
Then we have an obvious map $\wt{P}_\idot \to \overline{P}_\idot
\otimes I^\hdot$, and one easily shows that the induced map
$$
\vC'_\idot(G,\nu,E_\idot) \to \vC_\idot(G,\nu,E_\idot)
$$
is a quasiisomorphism. Thus the complex \eqref{c.dash.tate} also
computes Tate cohomology groups $\vH_\idot(G,E_\idot)$ (in fact,
historically, this was their original definition).

Now let us say that resolution data $\langle P_\idot,I^\hdot
\rangle$ are {\em multiplicative} if both $I^\hdot$ and
$\wt{P}_\idot$ are equipped with a structure of a unital associative
DG algebra. Note that multiplicative resolution data do exist for
any finite group $G$. Indeed, for $I^\hdot$, this is well-known, and
for $P_\idot$, it suffices to take a free $\Z[G]$-module $P_0$
equipped with a surjective map $d:P_0 \to \Z$, and let $P_i =
P_0^{\otimes i}$, $i \geq 1$, with the differential induced by $d$.

\begin{lemma}\label{dg.mult.cycl.le}
Assume given a set of multiplicative resolution data $\nu$ for the
cyclic group $\Z/p\Z$, and for any complex $V_\idot \in
C^{pf}_\idot(R)$ with tensor power $E_\idot = V_\idot^{\otimes_R
  p}$, let
$$
\wt{C}_\idot(\nu,V_\idot) = \tau^F_{[0,1]}\vC_\idot(\Z/p\Z,\nu,E_\idot),
$$
where the filtration $F^\hdot$ on $E_\idot$ is the same as in
Lemma~\ref{dg.tr.le}. Then the admissible functor
$\wt{C}_\idot(\nu,-):C^{pf}_\idot(R) \to C_\idot(R)$ has a natural
multiplicative structure in the sense of
Definition~\ref{DG.mult.def}.
\end{lemma}

\proof{} The tensor power $V^{\otimes_R p}$ has a trivial
multiplicative structure, and since $\nu$ is multiplicative, the
functor $\vC_\idot(\Z/p\Z,\nu,-)$ is also multiplicative. It
remains to notice that this multiplicativity is compatible with the
filtrations, and the filtered truncation functor
$\tau^F_{[0,1]}:CF_\idot(R) \to CF_\idot(R)$ also has an obvious
multiplicative structure.
\endproof

As an immediately corollary of Lemma~\ref{dg.mult.cycl.le}, we see
that the functor $C_\idot(\nu,V_\idot)$ obtained from
$\wt{C}_\idot(\nu,V_\idot)$ by the pullback square \eqref{sq.c.c}
is also multiplicative, so that we obtain a multiplicative DG
elementary extension $\phi_\idot^{\nu}$ of $R$ by $R^{(1)}$. Moreover, while the
periodic Tate complex $\vC_\idot(\Z/p\Z,E_\idot)$ is not of the form
\eqref{c.phi.tate} for any choice of the resolution data, it is of
the form \eqref{c.dash.tate} --- as resolution data, one takes the
pair of the standard periodic projective and injective resolutions
of $\Z$. Therefore we have a natural functorial quasiisomorphism
$$
\C_\idot(V_\idot) \to \C_\idot(\nu,V_\idot),
$$
where $C_\idot(V_\idot)$ is the functorial complex
\eqref{dg.c.v}. Taking restictions, we also obtain a map from the
elementary extension $\phi$ of \eqref{phi.df} to the restriction
$\rho(\phi_\idot^{\nu})$ of the DG elementary extension
$\phi_\idot^{\nu}$. An interested reader can easily check that this
map is multiplicative, so that $\phi^{\nu}_\idot$ indeed extends
$\phi$ to a multiplicative functor.

\subsection{Splittings.}

We will now show that under some assumptions, the situation for the
multiplicative DG elementary extension $\phi_\idot^{\nu}$ provided
by Lemma~\ref{dg.mult.cycl.le} is better than for a general
multiplicative DG elementary extension --- namely, we do have a good
multiplicative strict left DG splitting of the induced extension
$q_*(\phi_\idot^{\nu})$.

First of all, assume that the commutative ring $R$ is perfect ---
that is, the Frobenius map $\Fr:R \to R^{(1)}$ is bijective. In this
case, we have $\wt{C}_\idot(V) \cong C_\idot(V)$ and
$\wt{C}_\idot(V_\idot) \cong C_\idot(V_\idot)$, without the need to
apply the pullback square \eqref{sq.c.c}. Moreover, the augmentation
ideal in the second Witt vectors ring $W_2(R)=R'$ is generated by
$p$, so that we have $R \cong R'/p$.

Next, consider the category $C^{pf}_\idot(R')$, and for any complex
$V_\idot \in C^{pf}_\idot(R')$, denote
\begin{equation}\label{c.v.t}
C_\idot(V_\idot) =
\tau^F_{[0,1]}\vC_\idot(\Z/p\Z,V_\idot^{\otimes_{R'} p}),
\end{equation}
where as in Lemma~\ref{dg.tr.le}, the filtration $F^\hdot$ on
$V_\idot^{\otimes_{R'} p}$ is the $p$-th rescaling of the stupid
filtration.

\begin{lemma}\label{lft.ten.le}
Assume given a complex $V_\idot \in C^{pf}_\idot(R')$, let
$q:V_\idot \to V_\idot/p$ be the quotient map, and let $p:V_\idot/p
\to V_\idot$ be the embedding induced by the multiplication by $p$.
\begin{enumerate}
\item We have a natural exact sequence
\begin{equation}\label{p.2}
\begin{CD}
0 @>>> C_\idot(V_\idot/p) @>{p}>> C_\idot(V_\idot) @>{q}>>
C_\idot(V_\idot/p) @>>> 0
\end{CD}
\end{equation}
of complexes of $R'$-modules.
\item Let $q_i:\tau^F_{[i,i]}C_\idot(V_\idot) \to
  \tau^F_{[i,i]}C_\idot(V_\idot/p)$ be the map induced by the
  quotient map $q$ for $i=0,1$. Then $q_0$ is an isomorphism, and
  $q_1=0$.
\end{enumerate}
\end{lemma}

\proof{} For \thetag{i}, note that it suffices to prove that the
sequence becomes exact after taking the associated graded quotient
$\gr^i_F$ for an arbitrary integer $i$. We always have an exact
sequence
$$
\begin{CD}
0 @>>> \vC_\idot(\Z/p\Z,(V_\idot/p)^{\otimes_R p}) @>{p}>>
\vC_\idot(\Z/p\Z,V_\idot^{\otimes_{R'} p}) @>{q}>>\\
@.@>{q}>> \vC_\idot(\Z/p\Z,(V_\idot/p)^{\otimes_R p}) @>>> 0
\end{CD}
$$
of Tate complexes, so that by \eqref{gr.tau} and
Lemma~\ref{triv.le}, it suffices to check that the connecting
differential
$$
\delta_{i,j}:\vH_j(\Z/p\Z,\gr^i_F(V_\idot/p)^{\otimes_R p}) \to
\vH_{j-1}(\Z/p\Z,\gr^i_F(V_\idot/p)^{\otimes_R p})
$$
in the corresponding long exact sequence of homology vanishes for
any $i$ and $j=i,i+2$. Since the Tate complex is $2$-periodic, it
suffices to consider the case $j=i$. By Lemma~\ref{dg.tr.le}, we
have a functorial identification
\begin{equation}\label{iso}
\vH_i(\Z/p\Z,\gr^i_F(V_\idot/p)^{\otimes_R p}) \cong
\vH_{i-1}(\Z/p\Z,\gr^i_F(V_\idot/p)^{\otimes_R p}) \cong V_i^{(1)}/p,
\end{equation}
so that $\delta_{i,i}$ is an endomorphism of $V_i^{(1)}/p$. Then by
functoriality, we have replace the complex $V_\idot$ with $V_i[i]$,
so it suffices to consider the case when $V_\idot$ is concentrated
in a single homological degree.

For \thetag{ii}, note that since by Lemma~\ref{dg.tr.le}, the
filtration $F^\hdot$ induces a shift of the stupid filtration on
$\tau^F_{[j,j]}C_\idot(V_\idot)$ for any $j$, it again suffices to
prove the claim after passing to $\gr^i_F$. Moreover, if we know
that $\delta_{i,i}=0$, then both claims reduce to checking that
$\delta_{i,i+1}$ is an isomorphism, and by virtue of \eqref{iso}, it
suffices to check this for complexes concentrated in a single
homological degree.

Moreover, since both homology groups in \eqref{iso} are additive
with respect to $V_\idot$, we may further assume that $V_i \cong
R'$, so that $V_\idot = R'[i]$ for some integer $i$. Then
$V_\idot^{\otimes_{R'} p} \cong R'[ip]$, and the permutation
$\sigma$ acts by $(-1)^{i(p-1)}\id$. If $i(p-1)$ is even, then,
since the Tate complex \eqref{perio} is $2$-periodic, we have
\begin{equation}\label{shift}
\vC_\idot(\Z/p\Z,(R'[i])^{\otimes_{R'} p}) \cong
\vC_\idot(\Z/p\Z,R')[i].
\end{equation}
If $i(p-1)$ is odd, then necessarily $p=2$, and the differentials in
\eqref{perio} are given by $\id - \sigma$, $\id + \sigma$. Therefore
replacing $\sigma$ by $-\sigma$ is equivalent to shifting the
complex by $1$, and we still have the identification
\eqref{shift}. It shows that we may further assume that $i=0$. Then
the complex \eqref{perio} is the complex
$$
\begin{CD}
@>{0}>> R' @>{p}>> R' @>{0}>> R' @>{p}>> R' @>{0}>>,
\end{CD}
$$
with $d:\vC_i(\Z/p\Z,R') \to \vC_{i-1}(\Z/p\Z,R')$ given by $0$ for
even $i$ and $p$ for odd $i$, and both claims \thetag{i},
\thetag{ii} are obvious.
\endproof

We now note that the proof of Lemma~\ref{lft.ten.le} only depends on
the homology groups $\vH_\idot(\Z/p\Z,-)$, so it remains valid if we
change \eqref{c.v.t} by considering one of the different versions of
the Tate complex given in Subsection~\ref{mult.cycl.subs}. So, fix
once and for all a set $\nu$ of multiplicative resolution data for
the cyclic group $\Z/p\Z$, and let $\vC_\idot(\Z/p\Z,\nu,-) =
\vC_\idot(\Z/p\Z,-)$, with \eqref{c.v.t} reinterpreted accordingly
(we will not need the periodic complex \eqref{perio} anymore, so
there is no danger of confusion).

Now, for any complex $W_\idot \in C^{pf}_\idot(R)$, denote by
$\overline{C}^r_\idot(W_\idot) \subset C_\idot(W_\idot)$ the kernel
of the map $a:C_\idot(W_\idot) \to W_\idot$, and denote by
$\overline{C}^l_\idot = C_\idot(W_\idot)/b(W_\idot[1])$ the cokernel
of map $b:W_\idot[1] \to C_\idot(W_\idot)$. Then by
Lemma~\ref{lft.ten.le}~\thetag{i}, for any $V_\idot \in
C^{pf}_\idot(R')$, we can define complexes $C^r_\idot(V_\idot)$,
$C^l_\idot(V_\idot)$ by exact sequences
\begin{equation}\label{c.l.r.cycl}
\begin{CD}
0 @>>> \overline{C}^r_\idot(V_\idot/p) @>{p}>> C_\idot(V_\idot) @>>>
C^l_\idot(V_\idot) @>>> 0,\\
0 @>>> C^r_\idot(V_\idot) @>>> C_\idot(V_\idot) @>{q}>>
\overline{C}^l_\idot(V_\idot/p) @>>> 0.
\end{CD}
\end{equation}
By Lemma~\ref{lft.ten.le}~\thetag{ii}, $C^l_\idot(V_\idot)$
resp.\ $C^r_\idot(V_\idot)$ is a strict left resp.\ right DG
splitting of the DG elementary extension $C_\idot(V_\idot/p)$. Both
$C^r_\idot(V_\idot)$ and $C^l_\idot(V_\idot)$ are functorial in
$V_\idot$. The functor $C^l_\idot(-)$ has an obvious multiplicative
structure induced by the multiplicative structure on $C_\idot(-)$,
and the action of $C_\idot(-)$ on $C^r_\idot(-) \subset C_\idot(-)$
factors through its quotient $C^l_\idot(-)$ --- for any two
complexes $V_\idot,V'_\idot \in C^{pf}_\idot(R')$, we have
functorial action maps
\begin{equation}\label{alg.mod}
\begin{aligned}
C^l_\idot(V'_\idot) \otimes_{R'} C^l_\idot(V_\idot) &\to
C^l_\idot(V'_\idot \otimes_{R'} V_\idot),\\
C^l_\idot(V'_\idot) \otimes_{R'} C^r_\idot(V_\idot) &\to
C^r_\idot(V'_\idot \otimes_{R'} V_\idot).
\end{aligned}
\end{equation}
We will need one somewhat technical result on the functor
$C^l_\idot(-)$. Take the unity element $t_0=1 \in C^l_0(R')$, and as
in the proof of Lemma~\ref{dg.sec.le}, lift it to an element $s_0
\in C^l_0(\Cyl(R')[-1])$ such that $C^l_\idot(\alpha)(s_0)=t_0$ and
$ds_0=C^l_\idot(\beta)(t)$ for some $t_1 \in C_1(R'[-1])$. Then by
adjunction, $s_0$ induces a map
\begin{equation}\label{wt.s.df}
\wt{s}:V_0 \to C^l_0(V_\idot)
\end{equation}
for any complex $V_\idot \in C^{pf}_\idot(R')$, and this map is
functorial in $V_\idot$. Moreover, since $C^l_\idot(\alpha)(s_0)=1$,
we have
\begin{equation}\label{wt.s.v}
\wt{s}(v) = v^{\otimes p}
\end{equation}
for any closed $v \in V_0$, $dv=0$. If we project
$C^l_\idot(V_\idot)$ to $C_\idot(V_\idot/p)$, then the composition
map $V_0 \to C_0(V_\idot/p)$ factors through a functorial map
$$
\overline{s}:V_0/p \to C_0(V_\idot/p)
$$
corresponding to the image of the element $s_0$ in $C_0(R)$, and we
still have
\begin{equation}\label{ol.s.v}
\overline{s}(v) = v^{\otimes p}
\end{equation}
for any closed $v \in V_0/p$, $dv=0$. In particular, both $\wt{s}$
and $\overline{s}$ are multiplicative on closed elements. On the
whole $V_0$, the map $\wt{s}$ does not have to multiplicative, but
the following is sufficient for our purposes.

\begin{lemma}\label{fact.le}
\begin{enumerate}
\item Assume given two complexes $V_\idot,V'_\idot \in
  C^{pf}_\idot(R')$ and elements $v \in V_0$, $v' \in V'_0$ such
  that $dv=dv'=0 \mod p$. Then we have
\begin{equation}\label{s.mult.l}
\wt{s}(v \cdot v') = \wt{s}(v) \cdot \wt{s}(v') \in C^l_0(V_\idot
\otimes_{R'} V'_\idot).
\end{equation}
\item Assume in addition that $V_\idot = V'_\idot$ and $v=v' \mod
  p$. Then we have $\wt{s}(v) = \wt{s}(v')$.
\end{enumerate}
\end{lemma}

\proof{} For \thetag{i}, note that we can always replace the
complexes $V_\idot$, $V'_\idot$ with the $0$-th terms of their
stupid filtration; assume therefore that both lie in $C^{pf}_{\leq
  0}(R')$. Then so does the product $V''_\idot = V_\idot
\otimes_{R'} V'_\idot$, and the natural projection
$\lambda'':V''_\idot \to V''_0 \cong V_0 \otimes_{R'} V'_0$ is the
tensor product of the projections $\lambda:V_\idot \to V_0$,
$\lambda':V'_\idot \to V'_0$. Since $C^l_\idot(V''_\idot)$ is a
strict left DG splitting of $C_\idot(V''_\idot/p)$, the map
\begin{equation}\label{inj}
\begin{CD}
C^l_\idot(V''_\idot) @>{q \oplus C^l_\idot(\lambda'')}>>
C_\idot(V''_\idot/p) \oplus C^l_\idot(V''_0)
\end{CD}
\end{equation}
is injective in degree $0$ for dimension reasons. Therefore is
suffices to check \eqref{s.mult.l} after projecting to
$C_0(V''_\idot/p)$ and to $C^l_0(V''_0)$. For the former, note that
by assumption, $v$ and $v'$ are closed modulo $p$, so that
\eqref{s.mult.l} immediately follows from \eqref{ol.s.v}. For the
latter, note that then we can replace $V_\idot$, $V'_\idot$ with
$V_0$, $V'_0$, and then $v$ and $v'$ become closed, and
\eqref{s.mult.l} follows from \eqref{wt.s.v}.

For \thetag{ii}, note that it suffices to prove that
$\wt{s}(v)=\wt{s}(v')$ in the universal situation $V_\idot =
\Cyl(R'v_1 \oplus R'v_2)[-1]$, $v=v_1$, $v'=v_1 + pv_2$. In this
situation, the map \eqref{inj} for the complex $V_\idot$ is still
injective, and moreover, since $v=v' \mod p$, we have $q(\wt{s}(v))
= \overline{s}(v)=\overline{s}(v')=q(\wt{s}(v'))$, so that we may
further project to $V_0 = R'v_1 \oplus R'v_2$. Then \eqref{coc.eq.4}
and \eqref{c.exp} show that
\begin{equation}\label{vnsh}
\wt{s}(v')-\wt{s}(v) = \sum_{1 \leq i \leq p-1}p^i\sum_{s \in
  \overline{S}_i/G}\delta(\mu_{\kappa(s)}(v_1,v_2)),
\end{equation}
where $\overline{S}_i \subset \overline{S}$ is the subset
parametrizing monomials of degree $i$ in $v_2$ (and $p-i$ in
$v_1$). Since for any $i \geq 1$, $s \in \overline{S}_i/G$, $p^i$ is
divisible by $p$, while $\delta(\mu_{\kappa(s)}(v_1,v_2))$ lies in
$\overline{C}^r_0(V_0/p) \subset C_0(V_0/p)$, the
right-hand of \eqref{vnsh} gives $0$ in the quotient
$C^l_0(V_0)=C_0(V_0)/p\overline{C}^r_0(V_0/p)$.
\endproof

\subsection{Liftings.}

We can now prove the following surprising general property of the
functors $C^l_\idot$, $C^r_\idot$ of \eqref{c.l.r.cycl} (we note
that the idea for this is essentially due to V. Vologodsky). We need
to assume further that the perfect ring $R$ is a field.

\begin{prop}\label{lft.ten.prop}
The functorial complexes $C^l_\idot(V_\idot)$, $C^r_\idot(V_\idot)$
of \eqref{c.l.r.cycl} only depend on the quotient $V_\idot/p$, so
that the corresponding admissible functors
$$
C^l_\idot,C^r_\idot:C^{pf}_\idot(R') \to C_\idot(R')
$$
factor through the projection $q^*:C^{pf}_\idot(R') \to
C^{pf}_\idot(R)$, $q^*(V_\idot)=V_\idot/p$.
\end{prop}

\proof{} By definition, morphisms from $V_\idot$ to $V'_\idot$ in
the category $C^{pf}_\idot(R')$ are degree-$0$ classes
$$
f \in \Hom^0(V_\idot,V'_\idot)
$$
in the complex $\Hom^\hdot(V_\idot,V'_\idot)$ such that
$df=0$. Such a morphism acts on $V_\idot$ via the action map
\begin{equation}\label{act}
a:\Hom^\hdot(V_\idot,V'_\idot) \otimes_{R'} V_\idot \to V'_\idot.
\end{equation}
On the other hand, on $p$-th tensor powers, such a morphism $f$ acts
by $f^{\otimes p}$. Then \eqref{wt.s.v} shows that the morphisms
$C^l_\idot(f)$, $C^r_\idot(f)$ can be expressed in terms of the map
$\wt{s}$ of \eqref{wt.s.df} and the maps \eqref{alg.mod}. Namely,
for any $f:V_\idot \to V'_\idot$ and $c \in C^l_\idot(V_\idot)$, $c'
\in C^r_\idot(V_\idot)$, we have
\begin{equation}\label{f.acts}
C^l_\idot(f)(c) = C^l_\idot(a)(\wt{s}(f) \cdot c), \qquad
C^r_\idot(f)(c') = C^r_\idot(a)(\wt{s}(f) \cdot c'),
\end{equation}
where $\cdot$ denotes the product \eqref{alg.mod}, and $a$ is the
action map \eqref{act}.

Now, since $R$ is by assumption a field, the projection
$q^*:C^{pf}_\idot(R') \to C^{pf}_\idot(R)$ is essentially
surjective. Therefore one can describe the category
$C^{pf}_\idot(R)$ in the following way: objects are complexes
$V_\idot \in C^{pf}_\idot(R')$, morphisms from $V_\idot$ to
$V'_\idot$ are closed degree-$0$ classes
\begin{equation}\label{f.p}
f \in \Hom^0(V_\idot/p,V'_\idot/p)=\Hom^0(V_\idot,V'_\idot)/p
\end{equation}
in the complex $\Hom^\hdot(V_\idot,V'_\idot)/p$.  But by
Lemma~\ref{fact.le}, for any $V_\idot,V'_\idot \in
C^{pf}_\idot(R')$, the map $\wt{s}$ actually factors as
$$
\begin{CD}
\Hom^0(V_\idot,V'_\idot) @>{q}>>
\Hom^0(V_\idot,V'_\idot)/p @>{\overline{s}}>>
C^l_0(\Hom^\hdot(V_\idot,V'_\idot)),
\end{CD}
$$
and the map $\overline{s}$ is multiplicative on closed classes. Then
to factor the functors $C^l_\idot(-)$, $C^r_\idot(-)$ through $q^*$,
just let morphisms $f$ of \eqref{f.p} act by \eqref{f.acts}.
\endproof

\begin{remark}
The only place where we have used the assumption that $R$ is field
is in concluding that $q^*:C^{pf}_\idot(R') \to C^{pf}_\idot(R)$ is
essentially surjective. It might be that homological properties of
perfect rings are good enough to insure this in a more general
situation, but I have not pursued this.
\end{remark}

Now, by virtue of Proposition~\ref{lft.ten.prop}, we can redefine
the functors $C^l_\idot(-)$, $C^r_\idot(-)$ as admissible functors
$$
\C^l_\idot,C^r_\idot \in \B_\idot(R,R')
$$
from $C^{pf}_\idot(R)$ to $C_\idot(R')$, so that what we denoted
earlier by $\C^l_\idot(V_\idot)$, $\C^r_\idot(V_\idot)$ will now
become $\C^l_\idot(V_\idot/p)$, $\C^r_\idot(V_\idot/p)$. Then
$C^l_\idot$, $C^r_\idot$ are strict left resp.\ right DG splittings
of the extension $q_*C_\idot$, where $C_\idot$ is the cyclic powers
extension \eqref{dg.c.v} of Lemma~\ref{cycl.dg.le}. By the unicity
clause of Proposition~\ref{dg.Z.spl.prop}, the splitting $C^r_\idot$
must coincide with the one constructed in that Proposition in the
context of a general DG elementary extension. The splitting
$C^l_\idot$ is new: we do not know how to construct it in the
general case. However, for the cyclic powers extension $C_\idot$, it
is perfectly well-defined. Moreover, by construction, it has a
natural multiplicative structure. As in
Subsection~\ref{lft.comp.subs}, we can further extend to a functor
$$
C^l_\idot:C^{fl}_\idot(R) \to C_\idot(R'),
$$
and the extended functor inherits a multiplicative structure. In
particular, for any DG algebra $A_\idot$ over $R$, we have a
functorial map
\begin{equation}\label{c.rc}
C^l_\idot(A_\idot) \to q_*C_\idot(A_\idot)
\end{equation}
of DG algebras over $R'$. This opens the way to various DG versions
of Proposition~\ref{lft.mult.prop}. We will only record the
simplest possible result in this direction.

\begin{prop}
Assume given a DG algebra $A_\idot$ over $R$, and consider the DG
algebra $C_\idot(A_\idot)$ with the augmentation map
$a:C_\idot(A_\idot) \to A_\idot$. Assume further that there exists a
DG algebra $A'_\idot$ flat over $R'$ such that $A_\idot \cong
A'_\idot \otimes_{R'} R = A'_\idot/p$. Then there exists a DG
algebra $\wt{C}^l_\idot(A_\idot)$ over $R$ and a DG algebra map
$l:\wt{C}^l_\idot(A_\idot) \to C_\idot(A_\idot)$ such that the
composition map
$$
\begin{CD}
\wt{C}^l_\idot(A_\idot) @>{l}>> C_\idot(A_\idot) @>{a}>> A_\idot
\end{CD}
$$
is a quasiisomorphism.
\end{prop}

\proof{} We have the DG algebra $C^l_\idot(A_\idot)$ over $R'$ and
the DG algebra map \eqref{c.rc}. Moreover, $C^l_\idot(A_\idot)$ is a
multiplicative strict left DG splitting $\psi \in
\spl_\idot^l(\phi)$ of the multiplicative elementary extension
$q_*\phi_\idot=q_*C_\idot(A_\idot)$ in the category
$C_\idot(R')$. On the other hand, the DG algebra $A'_\idot$ is an
extension of $A_\idot$ by $A_\idot^{(1)}$, thus defines an object
$\eps$ in the extension groupoid
$\ex_\idot(A_\idot,A^{(1)}_\idot)$. Applying the difference functor
\eqref{dg.sum.fu.1}, we obtain an object
$$
\psi - \eps = \langle \wt{C}^l_\idot(A_\idot),l,b^l \rangle \in
\spl^l_\idot(q_*\phi).
$$
As in Proposition~\ref{lft.spl.prop}, $\wt{C}^l_\idot(A_\idot)$ is
actually a complex over $R=R'/p$, and as in
Proposition~\ref{lft.mult.prop}, it has a natural DG algebra
structure.
\endproof

\medskip

{\footnotesize
\noindent
{\sc
Steklov Math Institute, Algebraic Geometry section\\
\mbox{}\hspace{30mm}and\\
Laboratory of Algebraic Geometry, NRU HSE\\
\mbox{}\hspace{30mm}and\\
Center for Geometry and Physics, IBS, Pohang, Rep. of Korea
}

\smallskip

\noindent
{\em E-mail address\/}: {\tt kaledin@mi.ras.ru}
}
\end{document}